\documentclass[a4paper,final]{article}
\usepackage[dvips]{graphicx}
\usepackage{amsmath} % Extra environments for multiline displayed eq
\usepackage{amssymb} % Defines symbol names 
\makeatletter
\def\@maketitle{
  \vbox to 1.5in{
   \hsize\textwidth
    \linewidth\hsize
    \vspace*{1.5cm}
    \centering
    {\bfseries\huge \@title \par}
    \vskip 2em
    {\large \begin{tabular}[t]{c}\@author \end{tabular}\par}
    \vfill}    \vspace*{1.0cm}
}
\renewcommand\section{\@startsection {section}{1}{\z@}
     {.7\baselineskip plus\baselineskip}{.5\baselineskip}
                                   {\normalfont\Large\bfseries}}
\renewcommand\section{\@startsection {section}{1}{\z@}
      {.5\baselineskip\@plus.7\baselineskip}{.3\baselineskip}
                                   {\normalfont\Large\bfseries}}
\renewcommand\subsection{\@startsection{subsection}{2}{\z@}
       {.5\baselineskip\@plus.7\baselineskip}{.3\baselineskip}
                                   {\normalfont\large\bfseries}}
\renewcommand\subsubsection{\@startsection{subsubsection}{3}{\z@}
      {.5\baselineskip\@plus.7\baselineskip}{.3\baselineskip}
                                     {\normalfont\normalsize\bfseries}}
\makeatother

\renewenvironment{abstract}
  {\normalfont
    \list{}{\labelwidth0pt
      \leftmargin0pt \rightmargin\leftmargin
      \listparindent\parindent \itemindent0pt
      \parsep0pt
      
    }
    \item[\hskip\labelsep\bfseries\abstractname\enspace --] \itshape
}{
  \endlist}

\newcommand{\keywordsname}{Keywords}
\newenvironment{keywords}
  {\normalfont
    \list{}{\labelwidth0pt
      \leftmargin0pt \rightmargin\leftmargin
      \listparindent\parindent \itemindent0pt
      \parsep0pt
      }
    \item[\hskip\labelsep\bfseries\keywordsname:]}{\endlist}

\begin{document}

\pagestyle{myheadings}
%\markboth{Advances and Applications of DSmT  - F. Smarandache \& J. Dezert,  Editors - June 2003}{Advances and Applications of DSmT  - F. Smarandache \& J. Dezert, Editors - Oct. 7th,  2003}
\markboth{}{}

\title{Combining uncertain and paradoxical \\ evidences for DSm hybrid models\\
}
\author{Jean Dezert\\
ONERA\\
29 Avenue de la  Division Leclerc \\
92320 Ch\^{a}tillon, France.\\
Jean.Dezert@onera.fr\\
\and
Florentin Smarandache\\
Department of Mathematics\\
University of New Mexico\\
Gallup, NM 87301, U.S.A.\\
smarand@unm.edu}

\date{}

\maketitle
\vspace{2cm}

\begin{abstract}
This paper presents a general method for combining uncertain and paradoxical source of evidences for a wide class of fusion problems.
From the foundations of the Dezert-Smarandache Theory (DSmT) we show how the DSm rule of combination can be adapted to take into account all possible integrity constraints (if any) of the problem under consideration due to the true nature of elements/concepts involved into it. We show how the Shafer's model can be considered as a specific DSm hybrid model and be easily handled by our approach and a new efficient rule of combination different from the Dempster's rule is obtained. Several simple examples are also provided to show the efficiency and the generality of the approach proposed in this work. 
\end{abstract}



%********************
\section{Introduction}
%********************

A new theory of plausible and paradoxical reasoning (DSmT) has been
developed by the authors in the last two years in order to resolve problems that did not
work in Dempster-Shafer and other fusion theories.  According to each
model/problem of fusion occurring, we develop here a DSm hybrid rule
which combines two or more masses of independent sources of information and
takes care of restraints, i.e. of sets which might become empty at time $t_l$
or new sets/elements that might arise in the frame of discernment at time
$t_{l+1}$. DSm hybrid rule is applied in a real time when the hyper-power set  $D^\Theta$ changes (i.e. the set of all propositions built from elements of frame $\Theta$ with $\cup$ and $\cap$ operators - see\cite{Dezert_2003f} for details),
either increasing or decreasing its focal elements, or when even $\Theta$
decreases or increases influencing the $D^\Theta$ as well, thus the dynamicity
of our DSmT.\\

The paper introduces the reader to the independence of sources of
evidences, which needs to be deeper studied in the future, then defines
the models and the DSm hybrid rule, which is different from other rules of
combination such as Dempster's, Yager's, Smets', Dubois-Prade's and gives
seven numerical examples of applying the DSm hybrid rule in various models
and several examples of dynamicity of DSmT, then the Bayesian DSm hybrid
models mixture.

%************************************************
\section{On the independence of the sources of evidences}
%************************************************

The notion on independence of sources of evidences plays a major role in the development of efficient data fusion algorithms 
but is very difficult to formally establish when manipulating uncertain and paradoxical information. Some attempts to define the independence of uncertain sources of evidences have been proposed by P. Smets and al. in the Dempster-Shafer Theory (DST) and Transferable Belief Model in \cite{Yaghlane_Smets_1999,Yaghlane_Smets_Mellouli_2002a,Yaghlane_Smets_Mellouli_2002b} and by other authors in possibility theory \cite{Dawid_1998,Dawid_1999,Fonck_1994,Shenoy_1994,Studeny_1993}. In the following we consider that $n$ sources of evidences are independent if the internal mechanism by which each source provides its own basic belief assignment doesn't depend on the mechanisms of other sources (i.e. there is no internal relationship between all mechanisms) or if the sources don't share (even partially) same knowledge/experience to establish their own basic belief assignment. This definition doesn't exclude the possibility for independent sources to provide the same (numerical) basic belief assignments. The fusion of dependent uncertain and paradoxical sources is much more complicated because, one has first to identify precisely the piece of redundant information between sources in order to remove it before applying fusion rules. The problem of combination of dependent sources is under investigation.

%************************************************
\section{DSm rule of combination for free-DSm models}
%************************************************

\subsection{Definition of the free-DSm model $\mathcal{M}^f(\Theta)$}
%************************************************

Let consider a finite frame $\Theta=\{\theta_1,\ldots\theta_n\}$ of the fusion problem under consideration. We abandon the Shafer's model by assuming here that the fuzzy/vague/relative nature of elements $\theta_i$  $i=1,\ldots,n$ of $\Theta$ can be non-exclusive.
We assume also that  no refinement of $\Theta$ into a new  finer {\it{exclusive frame of discernment}}  $\Theta^{\text{ref}}$ is possible. This is the free-DSm model $\mathcal{M}^f(\Theta)$ which can be viewed as the opposite (if we don't introduce non-existential constraints - see next section) of Shafer's model, denoted $\mathcal{M}^0(\Theta)$ where all $\theta_i$ are forced to be exclusive and therefore fully discernable. 

\subsection{Example of a  free-DSm model}
%************************************************

Let consider the frame of the problem $\Theta=\{\theta_1,\theta_2,\theta_3\}$. The free Dedekind's lattice $D^\Theta=\{\alpha_0,\ldots,\alpha_{18}\}$ over $\Theta$ owns the following 19 elements \cite{Dezert_2003,Dezert_2003f}

\begin{equation*}
\begin{array}{|lll|}
\hline
& \text{Elements of } \; D^\Theta  \; \text{for} \; \mathcal{M}^f(\Theta)& \\
\hline
\alpha_0\triangleq \emptyset               & &                                 \\
\alpha_1\triangleq\theta_1\cap\theta_2\cap\theta_3  \neq \emptyset   & &    \\
\alpha_2\triangleq\theta_1\cap\theta_2 \neq \emptyset   & \alpha_3\triangleq\theta_1\cap\theta_3\neq \emptyset  &    \alpha_4\triangleq\theta_2\cap\theta_3\neq \emptyset                          \\
\alpha_5\triangleq(\theta_1\cup\theta_2)\cap\theta_3\neq \emptyset   & \alpha_6\triangleq(\theta_1\cup\theta_3)\cap\theta_2\neq \emptyset  &  \alpha_7\triangleq(\theta_2\cup\theta_3)\cap\theta_1\neq \emptyset  \\
  \alpha_8\triangleq\{(\theta_1\cap\theta_2)\cup\theta_3\} \cap(\theta_1\cup\theta_2)\neq \emptyset  & &  \\
\alpha_9\triangleq\theta_1\neq \emptyset       & \alpha_{10}\triangleq\theta_2\neq \emptyset   & \alpha_{11}\triangleq\theta_3\neq \emptyset                                                    \\
\alpha_{12}\triangleq(\theta_1\cap\theta_2)\cup\theta_3\neq \emptyset     & \alpha_{13}\triangleq(\theta_1\cap\theta_3)\cup\theta_2\neq \emptyset &  \alpha_{14}\triangleq(\theta_2\cap\theta_3)\cup\theta_1\neq \emptyset \\
\alpha_{15}\triangleq\theta_1\cup\theta_2\neq \emptyset     &      \alpha_{16}\triangleq\theta_1\cup\theta_3\neq \emptyset  &                 \alpha_{17}\triangleq\theta_2\cup\theta_3\neq \emptyset \\
\alpha_{18}\triangleq\theta_1\cup\theta_2\cup\theta_3\neq \emptyset    & &  \\
\hline
\end{array}
\end{equation*}

The free-DSm model $\mathcal{M}^f(\Theta)$ assumes that all elements $\alpha_i$, $i>0$, are non-empty. This corresponds to the following Venn diagram where in the Smarandache's codification "$i$" denotes the part of the diagram which belongs to $\theta_i$ only, "$ij$" denotes the part of the diagram which belongs to $\theta_i$ and $\theta_j$ only, "$ijk$" denotes the part of the diagram which belongs  to $\theta_i$ and $\theta_j$ and $\theta_k$ only, etc \cite{Dezert_2003f}. On such Venn diagram representation of the model, we emphasize the fact that all boundaries of intersections must be seen/interpreted as only vague boundaries just  because the nature of elements $\theta_i$ can be, in general, only  vague, relative and imprecise.
\begin{figure}[hbtp]
\centering
\includegraphics[width=4cm]{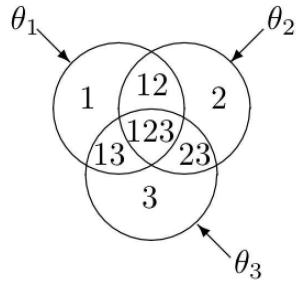}
\caption{Venn Diagram for $\mathcal{M}^f(\Theta)$ }
 \label{fig:2}
 \end{figure}

We now recall  the {\it{classical}} DSm rule of combination based $\mathcal{M}^f(\Theta)$ over the free Dedekind's lattice built from elements of $\Theta$ with $\cap$ and $\cup$ operators, i.e. the hyper-power set $D^\Theta$.

\subsection{Classical DSm rule for 2 sources for free-DSm models}
%***********************************************************************
For two {\it{independent }} uncertain and paradoxical sources of information (experts/bodies of evidence) providing generalized basic belief assignment $m_1(.)$ and $m_2(.)$ over $D^\Theta$ (or over any subset of $D^\Theta$), the classical DSm conjunctive rule of combination $m_{\mathcal{M}_f(\Theta)}(.)\triangleq [m_1\oplus m_2](.)$ is given by \cite{Dezert_2003}

\begin{equation}
\forall A\neq\emptyset\in D^\Theta, \quad m_{\mathcal{M}^f(\Theta)}(A)\triangleq [m_1\oplus m_2](A)=\sum_{\overset{X_1,X_2\in D^\Theta}{(X_1\cap X_2)=A}} m_1(X_1)m_2(X_2)
\label{eq:DSMClassic2}
\end{equation}

\noindent
$m_{\mathcal{M}^f(\Theta)}(\emptyset)=0$ by definition, unless otherwise specified in special cases when some source assigns a non-zero value to it (like in the Smets TBM approach \cite{Smets_1994}). This classic DSm rule of combination working on free-DSm models is commutative and associative. This rule, dealing with both uncertain and paradoxist/conflicting information, requires no normalization process and can always been applied.

\subsection{Classical DSm rule for $k\geq 2$ sources for free-DSm models}
%\subsection{Classical DSm rule of combination}
%************************************************

The above formula can be easily generalized for the free-DSm model $\mathcal{M}^f(\Theta)$ with $k\geq 2$ independent sources in the following way: 
%$\forall A\neq\emptyset \in D^\Theta$,
\begin{equation}
\forall A\neq\emptyset\in D^\Theta, \quad m_{\mathcal{M}^f(\Theta)}(A)\triangleq [m_1\oplus \ldots m_k](A)=
\sum_{\overset{X_1,\ldots,X_k\in D^\Theta}{(X_1\cap \ldots\cap X_k)=A}} \prod_{i=1}^{k}m_i(X_i)
\label{eq:DSMClassick}
\end{equation}
\noindent $m_{\mathcal{M}^f(\Theta)}(\emptyset)=0$ by definition, unless otherwise specified in special cases when some source assigns a non-zero value to it. This classic DSm rule of combination is still commutative and associative.

%************************************************
\section{Presentation of DSm hybrid models}
%************************************************

\subsection{Definition}
%************************************************

Let $\Theta$ be the general frame of the fusion problem under consideration with $n$ elements $\theta_1$, $\theta_2$, $\ldots$, $\theta_n$. A DSm hybrid model  $\mathcal{M}(\Theta)$ is defined from the free-DSm model $\mathcal{M}^f(\Theta)$ by introducing some {\it{integrity constraints}} on some elements $A$ of $D^\Theta$ if one knows with certainty the exact nature of the model corresponding to the problem under consideration. An integrity constraint on $A$ consists in forcing $A$ to be empty (vacuous element), and we will denote such constraint as $A\overset{\mathcal{M}}{\equiv}\emptyset$ which means that $A$ has been forced to $\emptyset$ through the model $\mathcal{M}(\Theta)$. 
This can be justified by the knowledge of the true nature of each element $\theta_i$ of $\Theta$. Indeed, in some fusion problems, some elements $\theta_i$ and $\theta_j$ of $\Theta$  can be fully discernable because they are truly exclusive while other elements cannot be refined into finer exclusive elements. Moreover, it is also possible that for some reason with some new knowledge on the problem, an element or several elements $\theta_i$ have to be forced to the empty set (specially if dynamical fusion problems are considered, i.e when $\Theta$ varies with space and time). For example, if we consider a list of three potential suspects into a police investigation, it can occur that, during the investigation, one of the suspect must be stray off  the initial frame of the problem if he can prove his gullibility with an ascertainable alibi. The initial basic belief masses provided by sources of information one had on the three suspects, must then be modified by taking into account this new knowledge on the model of the problem.\\

There exists several possible kinds of integrity constraints which can be introduced in any  free-DSm model $\mathcal{M}^f(\Theta)$ actually. The first kind of integrity constraint concerns {\it {exclusivity constraints}} by taking into account that some conjunctions of elements $\theta_i,\ldots , \theta_k$ are truly impossible (i.e. $\theta_i\cap\ldots\cap\theta_k\overset{\mathcal{M}}{\equiv}\emptyset$). The second kind of integrity constraint concerns the {\it{non-existential constraints}} by taking into account that some disjunctions of elements $\theta_i,\ldots,\theta_k$ are also truly impossible (i.e. 
$\theta_i\cup\ldots\cup\theta_k\overset{\mathcal{M}}{\equiv}\emptyset$). We exclude from our presentation the completely degenerate case corresponding to the constraint $\theta_1\cup\ldots\cup\theta_n\overset{\mathcal{M}}{\equiv}\emptyset$ (total ignorance) because there is no way and interest to treat such vacuous problem. In such degenerate case, we can just set $m(\emptyset)\triangleq 1$ which is useless because the problem remains vacuous and $D^\Theta$ reduces to $\emptyset$.
The last kind of possible integrity constraint is a mixture of the two previous ones, like for example $(\theta_i\cap\theta_j)\cup \theta_k$ or any other hybrid proposition/element of $D^\Theta$ involving both $\cap$ and $\cup$ operators such that at least one element $\theta_k$ is subset of the constrained proposition.
From any  $\mathcal{M}^f(\Theta)$, we can thus build several DSm hybrid models depending on the number of integrity constraints one wants to fully characterize the nature of the problem. The introduction of a given integrity constraint $A\overset{\mathcal{M}}{\equiv}\emptyset \in D^\Theta$ implies necessarily the set of inner constraints $B\overset{\mathcal{M}}{\equiv}\emptyset$ for all $B\subset A$. Moreover the introduction of two integrity constraints, say on $A$ and $B$ in $D^\Theta$ implies also necessarily the constraint on the emptiness of the disjunction $A\cup B$ which belongs also to $D^\Theta$ (because $D^\Theta$ is close under $\cap$ and $\cup$ operators). This implies the emptiness of all $C\in D^\Theta$ such that $C\subset (A\cup B)$. Same remark has to be extended for the case of the introduction of $n$ integrity constraints as well. The Shafer's model is the unique and most constrained DSm hybrid model including all possible exclusivity constraints {\it{without non-existential constraint  }} since all $\theta_i\neq\emptyset \in \Theta$ are forced to be mutually exclusive. The Shafer's model is denoted  $\mathcal{M}^{0}(\Theta)$ in the sequel. We denote by $\boldsymbol{\emptyset}_ {\mathcal{M}}$ the set of elements of $D^\Theta$ which have been forced to be empty in the DSm hybrid model $\mathcal{M}$.

\subsection{Example 1 : DSm hybrid model with an exclusivity constraint}
%**************************************************************************

Let $\Theta=\{\theta_1,\theta_2,\theta_3\}$ be the general frame of the problem under consideration and let consider the following DSm hybrid model $\mathcal{M}_1(\Theta)$ built by introducing the following exclusivity constraint $\alpha_1\triangleq\theta_1\cap\theta_2\cap\theta_3  \overset{\mathcal{M}_1}{\equiv}\emptyset$. This exclusivity constraint implies however no other constraint because $\alpha_1$ doesn't contain other elements of $D^\Theta$ but itself.
Therefore, one has now the following set of elements for $D^\Theta$
\begin{equation*}
\begin{array}{|lll|}
\hline
& \text{Elements of } \; D^\Theta  \; \text{for} \; \mathcal{M}_1(\Theta)& \\
\hline
\alpha_0\triangleq \emptyset               & &                                 \\
\alpha_1\triangleq\theta_1\cap\theta_2\cap\theta_3 \overset{\mathcal{M}_1}{\equiv}\emptyset   & &    \\
\alpha_2\triangleq\theta_1\cap\theta_2 \neq \emptyset   & \alpha_3\triangleq\theta_1\cap\theta_3\neq \emptyset  &    \alpha_4\triangleq\theta_2\cap\theta_3\neq \emptyset                          \\
\alpha_5\triangleq(\theta_1\cup\theta_2)\cap\theta_3\neq \emptyset   & \alpha_6\triangleq(\theta_1\cup\theta_3)\cap\theta_2\neq \emptyset  &  \alpha_7\triangleq(\theta_2\cup\theta_3)\cap\theta_1\neq \emptyset  \\
  \alpha_8\triangleq\{(\theta_1\cap\theta_2)\cup\theta_3\} \cap(\theta_1\cup\theta_2)\neq \emptyset  & &  \\
\alpha_9\triangleq\theta_1\neq \emptyset       & \alpha_{10}\triangleq\theta_2\neq \emptyset   & \alpha_{11}\triangleq\theta_3\neq \emptyset                                                    \\
\alpha_{12}\triangleq(\theta_1\cap\theta_2)\cup\theta_3\neq \emptyset     & \alpha_{13}\triangleq(\theta_1\cap\theta_3)\cup\theta_2\neq \emptyset &  \alpha_{14}\triangleq(\theta_2\cap\theta_3)\cup\theta_1\neq \emptyset \\
\alpha_{15}\triangleq\theta_1\cup\theta_2\neq \emptyset     &      \alpha_{16}\triangleq\theta_1\cup\theta_3\neq \emptyset  &                 \alpha_{17}\triangleq\theta_2\cup\theta_3\neq \emptyset \\
\alpha_{18}\triangleq\theta_1\cup\theta_2\cup\theta_3\neq \emptyset    & &  \\
\hline
\end{array}
\end{equation*}

Hence the initial basic belief mass over $D^\Theta$ has to be transferred over the new constrained hyper-power set $D^\Theta(\mathcal{M}_1(\Theta))$ with the 18 elements defined just above. The mechanism for the transfer of basic belief masses from $D^\Theta$ onto $D^\Theta(\mathcal{M}_1(\Theta))$ will be obtained by the DSm hybrid rule of combination presented in the sequel.

\subsection{Example 2 : DSm hybrid model with another exclusivity constraint}
%************************************************************************************

As second example for DSm hybrid model $\mathcal{M}_2(\Theta)$, let consider $\Theta=\{\theta_1,\theta_2,\theta_3\}$ and the following exclusivity constraint $\alpha_2\triangleq\theta_1\cap\theta_2 \overset{\mathcal{M}_2}{\equiv}\emptyset$. This constraint implies also $\alpha_1\triangleq\theta_1\cap\theta_2\cap\theta_3\overset{\mathcal{M}_2}{\equiv}\emptyset$ since $\alpha_1\subset\alpha_2$. Therefore, one has now the following set of elements for $D^\Theta (\mathcal{M}_2(\Theta))$
\begin{equation*}
\begin{array}{|lll|}
\hline
& \text{Elements of } \; D^\Theta  \; \text{for} \; \mathcal{M}_2(\Theta)& \\
\hline
\alpha_0\triangleq \emptyset               & &                                 \\
\alpha_1\triangleq\theta_1\cap\theta_2\cap\theta_3  \overset{\mathcal{M}_2}{\equiv}\emptyset   & &    \\
\alpha_2\triangleq\theta_1\cap\theta_2 \overset{\mathcal{M}_2}{\equiv}\emptyset   & \alpha_3\triangleq\theta_1\cap\theta_3\neq \emptyset  &    \alpha_4\triangleq\theta_2\cap\theta_3\neq \emptyset                          \\
\alpha_5\triangleq(\theta_1\cup\theta_2)\cap\theta_3\neq \emptyset   & \alpha_6\triangleq(\theta_1\cup\theta_3)\cap\theta_2\overset{\mathcal{M}_2}{\equiv} 
\alpha_4\neq \emptyset  &  \alpha_7\triangleq(\theta_2\cup\theta_3)\cap\theta_1\overset{\mathcal{M}_2}{\equiv}
\alpha_3\neq \emptyset  \\
  \alpha_8\triangleq\{(\theta_1\cap\theta_2)\cup\theta_3\} \cap(\theta_1\cup\theta_2) \overset{\mathcal{M}_2}{\equiv}\alpha_5 \neq \emptyset  & &  \\
\alpha_9\triangleq\theta_1\neq \emptyset       & \alpha_{10}\triangleq\theta_2\neq \emptyset   & \alpha_{11}\triangleq\theta_3\neq \emptyset                                                    \\
\alpha_{12}\triangleq(\theta_1\cap\theta_2)\cup\theta_3 \overset{\mathcal{M}_2}{\equiv} \alpha_{11}\neq \emptyset     & \alpha_{13}\triangleq(\theta_1\cap\theta_3)\cup\theta_2\neq \emptyset &  \alpha_{14}\triangleq(\theta_2\cap\theta_3)\cup\theta_1\neq \emptyset \\
\alpha_{15}\triangleq\theta_1\cup\theta_2\neq \emptyset     &      \alpha_{16}\triangleq\theta_1\cup\theta_3\neq \emptyset  &                 \alpha_{17}\triangleq\theta_2\cup\theta_3\neq \emptyset \\
\alpha_{18}\triangleq\theta_1\cup\theta_2\cup\theta_3\neq \emptyset    & &  \\
\hline
\end{array}
\end{equation*}
\noindent Note that in this case several non-empty elements of $D^\Theta(\mathcal{M}_2(\Theta))$ coincide because of the constraint ($\alpha_6 \overset{\mathcal{M}_2}{\equiv} \alpha_4$,
$\alpha_7 \overset{\mathcal{M}_2}{\equiv} \alpha_3$, $\alpha_8 \overset{\mathcal{M}_2}{\equiv} \alpha_5$, $\alpha_{12} \overset{\mathcal{M}_2}{\equiv} \alpha_{11}$). $D^\Theta(\mathcal{M}_2(\Theta))$ has now only 13 different elements. Note that the introduction of both constraints $\alpha_1\triangleq\theta_1\cap\theta_2\cap\theta_3  \overset{\mathcal{M}_2}{\equiv}\emptyset$ and $\alpha_2\triangleq\theta_1\cap\theta_2 \overset{\mathcal{M}_2}{\equiv}\emptyset$ doesn't change the construction of $D^\Theta(\mathcal{M}_2(\Theta))$ because $\alpha_1\subset\alpha_2$.

\clearpage
\newpage

\subsection{Example 3 : DSm hybrid model with another exclusivity constraint}
%************************************************************************************

As third example  for DSm hybrid model $\mathcal{M}_3(\Theta)$, let consider $\Theta=\{\theta_1,\theta_2,\theta_3\}$ and the following exclusivity constraint $\alpha_6\triangleq(\theta_1\cup\theta_3)\cap\theta_2\overset{\mathcal{M}_3}{\equiv}\emptyset$. This constraint implies now $\alpha_1\triangleq\theta_1\cap\theta_2\cap\theta_3\overset{\mathcal{M}_3}{\equiv}\emptyset$ since $\alpha_1\subset\alpha_6$, but also
$\alpha_2\triangleq\theta_1\cap\theta_2 \overset{\mathcal{M}_3}{\equiv}\emptyset$ because $\alpha_2\subset\alpha_6$ and $\alpha_4\triangleq\theta_2\cap\theta_3\overset{\mathcal{M}_3}{\equiv}\emptyset$ because $\alpha_4\subset\alpha_6$.
Therefore, one has now the following set of elements for $D^\Theta (\mathcal{M}_3(\Theta))$
\begin{equation*}
\begin{array}{|lll|}
\hline
& \text{Elements of } \; D^\Theta  \; \text{for} \; \mathcal{M}_3(\Theta)& \\
\hline
\alpha_0\triangleq \emptyset               & &                                 \\
\alpha_1\triangleq\theta_1\cap\theta_2\cap\theta_3 \overset{\mathcal{M}_3}{\equiv}\emptyset   & &    \\
\alpha_2\triangleq\theta_1\cap\theta_2 \overset{\mathcal{M}_3}{\equiv}\emptyset   & \alpha_3\triangleq\theta_1\cap\theta_3\neq \emptyset  &    \alpha_4\triangleq\theta_2\cap\theta_3\overset{\mathcal{M}_3}{\equiv}\emptyset                          \\
\alpha_5\triangleq(\theta_1\cup\theta_2)\cap\theta_3 \overset{\mathcal{M}_3}{\equiv}\alpha_{3} \neq \emptyset   & \alpha_6\triangleq(\theta_1\cup\theta_3)\cap\theta_2\overset{\mathcal{M}_3}{\equiv}\emptyset &  \alpha_7\triangleq(\theta_2\cup\theta_3)\cap\theta_1\overset{\mathcal{M}_3}{\equiv}\alpha_{3} \neq \emptyset  \\
  \alpha_8\triangleq\{(\theta_1\cap\theta_2)\cup\theta_3\} \cap(\theta_1\cup\theta_2)\overset{\mathcal{M}_3}{\equiv}\alpha_{5}\neq \emptyset  & &  \\
\alpha_9\triangleq\theta_1\neq \emptyset       & \alpha_{10}\triangleq\theta_2\neq \emptyset   & \alpha_{11}\triangleq\theta_3\neq \emptyset                                                    \\
\alpha_{12}\triangleq(\theta_1\cap\theta_2)\cup\theta_3\overset{\mathcal{M}_3}{\equiv}\alpha_{11}\neq \emptyset     & \alpha_{13}\triangleq(\theta_1\cap\theta_3)\cup\theta_2\neq \emptyset &  \alpha_{14}\triangleq(\theta_2\cap\theta_3)\cup\theta_1\overset{\mathcal{M}_3}{\equiv}\alpha_{9}\neq \emptyset \\
\alpha_{15}\triangleq\theta_1\cup\theta_2\neq \emptyset     &      \alpha_{16}\triangleq\theta_1\cup\theta_3\neq \emptyset  &                 \alpha_{17}\triangleq\theta_2\cup\theta_3\neq \emptyset \\
\alpha_{18}\triangleq\theta_1\cup\theta_2\cup\theta_3\neq \emptyset    & &  \\
\hline
\end{array}
\end{equation*}

\noindent
$D^\Theta(\mathcal{M}_3(\Theta))$ has now only 10 different elements.

\subsection{Example 4 : DSm hybrid model with all exclusivity constraints}
%*******************************************************************************

As fourth example for DSm hybrid model $\mathcal{M}_4(\Theta)$, let consider $\Theta=\{\theta_1,\theta_2,\theta_3\}$ and the following exclusivity constraint $\alpha_8\triangleq\{(\theta_1\cap\theta_2)\cup\theta_3\} \cap(\theta_1\cup\theta_2)\overset{\mathcal{M}_4}{\equiv}\emptyset$. 
This model corresponds actually to the Shafer's model $\mathcal{M}^0(\Theta)$ because this constraint includes all possible exclusivity constraints between elements $\theta_i$, $i=1,2,3$ since $\alpha_1\triangleq\theta_1\cap\theta_2\cap\theta_3 \subset \alpha_8$, $\alpha_2\triangleq\theta_1\cap\theta_2 \subset \alpha_8$, $\alpha_3\triangleq\theta_1\cap\theta_3 \subset \alpha_8$ and $\alpha_4\triangleq\theta_2\cap\theta_3 \subset \alpha_8$.
Therefore, one has now the following set of elements for $D^\Theta (\mathcal{M}_4(\Theta))$
\begin{equation*}
\begin{array}{|lll|}
\hline
& \text{Elements of } \; D^\Theta  \; \text{for} \; \mathcal{M}_4(\Theta)& \\
\hline
\alpha_0\triangleq \emptyset               & &                                 \\
\alpha_1\triangleq\theta_1\cap\theta_2\cap\theta_3 \overset{\mathcal{M}_4}{\equiv}\emptyset   & &    \\
\alpha_2\triangleq\theta_1\cap\theta_2 \overset{\mathcal{M}_4}{\equiv}\emptyset   & \alpha_3\triangleq\theta_1\cap\theta_3 \overset{\mathcal{M}_4}{\equiv}\emptyset  &    \alpha_4\triangleq\theta_2\cap\theta_3\overset{\mathcal{M}_4}{\equiv}\emptyset                          \\
\alpha_5\triangleq(\theta_1\cup\theta_2)\cap\theta_3\overset{\mathcal{M}_4}{\equiv}\emptyset   & \alpha_6\triangleq(\theta_1\cup\theta_3)\cap\theta_2\overset{\mathcal{M}_4}{\equiv}\emptyset&  \alpha_7\triangleq(\theta_2\cup\theta_3)\cap\theta_1\overset{\mathcal{M}_4}{\equiv}\emptyset  \\
  \alpha_8\triangleq\{(\theta_1\cap\theta_2)\cup\theta_3\} \cap(\theta_1\cup\theta_2)\overset{\mathcal{M}_4}{\equiv}\emptyset  & &  \\
\alpha_9\triangleq\theta_1\neq \emptyset       & \alpha_{10}\triangleq\theta_2\neq \emptyset   & \alpha_{11}\triangleq\theta_3\neq \emptyset                                                    \\
\alpha_{12}\triangleq(\theta_1\cap\theta_2)\cup\theta_3\overset{\mathcal{M}_4}{\equiv}\alpha_{11}\neq \emptyset     & \alpha_{13}\triangleq(\theta_1\cap\theta_3)\cup\theta_2\overset{\mathcal{M}_4}{\equiv}\alpha_{10}\neq \emptyset &  \alpha_{14}\triangleq(\theta_2\cap\theta_3)\cup\theta_1\overset{\mathcal{M}_4}{\equiv}\alpha_{9}\neq \emptyset \\
\alpha_{15}\triangleq\theta_1\cup\theta_2\neq \emptyset     &      \alpha_{16}\triangleq\theta_1\cup\theta_3\neq \emptyset  &                 \alpha_{17}\triangleq\theta_2\cup\theta_3\neq \emptyset \\
\alpha_{18}\triangleq\theta_1\cup\theta_2\cup\theta_3\neq \emptyset    & &  \\
\hline
\end{array}
\end{equation*}

$D^\Theta (\mathcal{M}_4(\Theta))$ has now $2^{\vert \Theta \vert}=8$ different elements and coincides obviously with the classical power set $2^\Theta$. This corresponds to the Shafer's model and serves as foundation for the Dempster-Shafer Theory.

\clearpage
\newpage

\subsection{Example 5 : DSm hybrid model with a non-existential constraint}
%*******************************************************************************

As fifth example for DSm hybrid model $\mathcal{M}_5(\Theta)$, let consider $\Theta=\{\theta_1,\theta_2,\theta_3\}$ and the following {\it{non-existential}} constraint $
\alpha_9\triangleq\theta_1\overset{\mathcal{M}_5}{\equiv}\emptyset$. 
In other words, we remove $\theta_1$ from the initial frame $\Theta=\{\theta_1,\theta_2,\theta_3\}$.
This non-existential constraint implies $\alpha_1\triangleq\theta_1\cap\theta_2\cap\theta_3 \overset{\mathcal{M}_5}{\equiv}\emptyset$, $\alpha_2\triangleq\theta_1\cap\theta_2 \overset{\mathcal{M}_5}{\equiv}\emptyset$, $\alpha_3\triangleq\theta_1\cap\theta_3 \overset{\mathcal{M}_5}{\equiv}\emptyset$ and $\alpha_7\triangleq(\theta_2\cup\theta_3)\cap\theta_1\overset{\mathcal{M}_5}{\equiv}\emptyset $.
Therefore, one has now the following set of elements for $D^\Theta (\mathcal{M}_5(\Theta))$
\begin{equation*}
\begin{array}{|lll|}
\hline
& \text{Elements of } \; D^\Theta  \; \text{for} \; \mathcal{M}_5(\Theta)& \\
\hline
\alpha_0\triangleq \emptyset               & &                                 \\
\alpha_1\triangleq\theta_1\cap\theta_2\cap\theta_3 \overset{\mathcal{M}_5}{\equiv}\emptyset   & &    \\
\alpha_2\triangleq\theta_1\cap\theta_2 \overset{\mathcal{M}_5}{\equiv}\emptyset   & \alpha_3\triangleq\theta_1\cap\theta_3 \overset{\mathcal{M}_5}{\equiv}\emptyset  &    \alpha_4\triangleq\theta_2\cap\theta_3\neq\emptyset                          \\
\alpha_5\triangleq(\theta_1\cup\theta_2)\cap\theta_3\overset{\mathcal{M}_5}{\equiv}\alpha_4\neq\emptyset   & \alpha_6\triangleq(\theta_1\cup\theta_3)\cap\theta_2\overset{\mathcal{M}_5}{\equiv}\alpha_4\neq\emptyset&  \alpha_7\triangleq(\theta_2\cup\theta_3)\cap\theta_1\overset{\mathcal{M}_5}{\equiv}\emptyset  \\
  \alpha_8\triangleq\{(\theta_1\cap\theta_2)\cup\theta_3\} \cap(\theta_1\cup\theta_2)\overset{\mathcal{M}_5}{\equiv}\alpha_4\neq\emptyset  & &  \\
\alpha_9\triangleq\theta_1 \overset{\mathcal{M}_5}{\equiv}\emptyset      & \alpha_{10}\triangleq\theta_2\neq \emptyset   & \alpha_{11}\triangleq\theta_3\neq \emptyset                                                    \\
\alpha_{12}\triangleq(\theta_1\cap\theta_2)\cup\theta_3\overset{\mathcal{M}_5}{\equiv}\alpha_{11}\neq \emptyset     & \alpha_{13}\triangleq(\theta_1\cap\theta_3)\cup\theta_2\overset{\mathcal{M}_5}{\equiv}\alpha_{10}\neq \emptyset &  \alpha_{14}\triangleq(\theta_2\cap\theta_3)\cup\theta_1\overset{\mathcal{M}_5}{\equiv}\alpha_{4}\neq \emptyset \\
\alpha_{15}\triangleq\theta_1\cup\theta_2\overset{\mathcal{M}_5}{\equiv}\alpha_{10}\neq \emptyset     &      \alpha_{16}\triangleq\theta_1\cup\theta_3\overset{\mathcal{M}_5}{\equiv}\alpha_{11}\neq \emptyset  &                 \alpha_{17}\triangleq\theta_2\cup\theta_3\neq \emptyset \\
\alpha_{18}\triangleq\theta_1\cup\theta_2\cup\theta_3\overset{\mathcal{M}_5}{\equiv}\alpha_{17}\neq \emptyset    & &  \\
\hline
\end{array}
\end{equation*}

$D^\Theta (\mathcal{M}_5(\Theta))$ has now 5 different elements and coincides obviously with the hyper-power set $D^{\Theta\setminus\theta_1}$.

\subsection{Example 6 : DSm hybrid model with two non-existential constraints}
%**************************************************************************************

As sixth example for DSm hybrid model $\mathcal{M}_6(\Theta)$, let consider $\Theta=\{\theta_1,\theta_2,\theta_3\}$ and the following two {\it{non-existential}} constraints $
\alpha_9\triangleq\theta_1\overset{\mathcal{M}_6}{\equiv}\emptyset$ and $
\alpha_{10}\triangleq\theta_2\overset{\mathcal{M}_6}{\equiv}\emptyset$. Actually, these two constraints are equivalent to choose only the following constraint $\alpha_{15}\triangleq\theta_1\cup\theta_2\overset{\mathcal{M}_5}{\equiv}\emptyset$. 
In other words, we remove now both $\theta_1$ and $\theta_2$ from the initial frame $\Theta=\{\theta_1,\theta_2,\theta_3\}$. These non-existential constraints implies now $\alpha_1\triangleq\theta_1\cap\theta_2\cap\theta_3 \overset{\mathcal{M}_6}{\equiv}\emptyset$, $\alpha_2\triangleq\theta_1\cap\theta_2 \overset{\mathcal{M}_6}{\equiv}\emptyset$, $\alpha_3\triangleq\theta_1\cap\theta_3 \overset{\mathcal{M}_6}{\equiv}\emptyset$,
$\alpha_4\triangleq\theta_2\cap\theta_3 \overset{\mathcal{M}_6}{\equiv}\emptyset$,
$\alpha_5\triangleq(\theta_1\cup\theta_2)\cap\theta_3\overset{\mathcal{M}_6}{\equiv}\emptyset$,
$\alpha_6\triangleq(\theta_1\cup\theta_3)\cap\theta_2\overset{\mathcal{M}_6}{\equiv}\emptyset$,
$\alpha_7\triangleq(\theta_2\cup\theta_3)\cap\theta_1\overset{\mathcal{M}_6}{\equiv}\emptyset $,
$ \alpha_8\triangleq\{(\theta_1\cap\theta_2)\cup\theta_3\} \cap(\theta_1\cup\theta_2)\overset{\mathcal{M}_6}{\equiv}\emptyset$,
$\alpha_{13}\triangleq(\theta_1\cap\theta_3)\cup\theta_2\overset{\mathcal{M}_6}{\equiv}\emptyset$,
$\alpha_{14}\triangleq(\theta_2\cap\theta_3)\cup\theta_1\overset{\mathcal{M}_6}{\equiv}\emptyset$
.
Therefore, one has now the following set of elements for $D^\Theta (\mathcal{M}_6(\Theta))$
\begin{equation*}
\begin{array}{|lll|}
\hline
& \text{Elements of } \; D^\Theta  \; \text{for} \; \mathcal{M}_6(\Theta)& \\
\hline
\alpha_0\triangleq \emptyset               & &                                 \\
\alpha_1\triangleq\theta_1\cap\theta_2\cap\theta_3 \overset{\mathcal{M}_6}{\equiv}\emptyset   & &    \\
\alpha_2\triangleq\theta_1\cap\theta_2 \overset{\mathcal{M}_6}{\equiv}\emptyset   & \alpha_3\triangleq\theta_1\cap\theta_3 \overset{\mathcal{M}_6}{\equiv}\emptyset  &    \alpha_4\triangleq\theta_2\cap\theta_3\overset{\mathcal{M}_6}{\equiv}\emptyset                          \\
\alpha_5\triangleq(\theta_1\cup\theta_2)\cap\theta_3\overset{\mathcal{M}_6}{\equiv}\emptyset   & \alpha_6\triangleq(\theta_1\cup\theta_3)\cap\theta_2\overset{\mathcal{M}_6}{\equiv}\emptyset &  \alpha_7\triangleq(\theta_2\cup\theta_3)\cap\theta_1\overset{\mathcal{M}_6}{\equiv}\emptyset  \\
  \alpha_8\triangleq\{(\theta_1\cap\theta_2)\cup\theta_3\} \cap(\theta_1\cup\theta_2)\overset{\mathcal{M}_6}{\equiv}\emptyset  & &  \\
\alpha_9\triangleq\theta_1 \overset{\mathcal{M}_6}{\equiv}\emptyset      & \alpha_{10}\triangleq\theta_2\overset{\mathcal{M}_6}{\equiv}\emptyset  & \alpha_{11}\triangleq\theta_3\neq \emptyset                                                    \\
\alpha_{12}\triangleq(\theta_1\cap\theta_2)\cup\theta_3\overset{\mathcal{M}_6}{\equiv}\alpha_{11}\neq \emptyset     & \alpha_{13}\triangleq(\theta_1\cap\theta_3)\cup\theta_2\overset{\mathcal{M}_6}{\equiv}\emptyset &  \alpha_{14}\triangleq(\theta_2\cap\theta_3)\cup\theta_1\overset{\mathcal{M}_6}{\equiv}\emptyset \\
\alpha_{15}\triangleq\theta_1\cup\theta_2\overset{\mathcal{M}_6}{\equiv}\emptyset    &      \alpha_{16}\triangleq\theta_1\cup\theta_3\overset{\mathcal{M}_6}{\equiv}\alpha_{11}\neq \emptyset  &                 \alpha_{17}\triangleq\theta_2\cup\theta_3\overset{\mathcal{M}_6}{\equiv}\alpha_{11}\neq \emptyset\\
\alpha_{18}\triangleq\theta_1\cup\theta_2\cup\theta_3\overset{\mathcal{M}_6}{\equiv}\alpha_{11}\neq \emptyset    & &  \\
\hline
\end{array}
\end{equation*}

$D^\Theta (\mathcal{M}_6(\Theta))$ reduces now to only two  different elements $\emptyset$ and $\theta_3$. $D^\Theta (\mathcal{M}_6(\Theta))$  coincides obviously with the hyper-power set $D^{\Theta\setminus\{\theta_1,\theta_2\}}$. Because there exists only one possible non empty element in $D^\Theta (\mathcal{M}_6(\Theta))$, such kind of problem is called a trivial problem. If one now introduces all non-existential constraints in free-DSm model, then the initial problem reduces to a {\it{vacuous problem}} also called  {\it{impossible problem}} corresponding to $m(\emptyset)\equiv1$ (no problem at all since the problem doesn't not exist now !!!). Such kinds of trivial or vacuous  {\it{problems}} are not considered anymore in the sequel since they present no real interest for engineering data fusion problems.

\subsection{Example 7 : DSm hybrid model with a mixed constraint}
%***********************************************************************

As seventh example for DSm hybrid model $\mathcal{M}_7(\Theta)$, let consider $\Theta=\{\theta_1,\theta_2,\theta_3\}$ and the following {\it{mixed exclusivity and  non-existential}} constraint 
$\alpha_{12}\triangleq(\theta_1\cap\theta_2)\cup\theta_3\overset{\mathcal{M}_7}{\equiv}\emptyset$. This mixed constraint implies
$\alpha_1\triangleq\theta_1\cap\theta_2\cap\theta_3 \overset{\mathcal{M}_7}{\equiv}\emptyset$,
$\alpha_2\triangleq\theta_1\cap\theta_2 \overset{\mathcal{M}_7}{\equiv}\emptyset$,
$\alpha_3\triangleq\theta_1\cap\theta_3 \overset{\mathcal{M}_7}{\equiv}\emptyset$,
$\alpha_4\triangleq\theta_2\cap\theta_3 \overset{\mathcal{M}_7}{\equiv}\emptyset$,
$\alpha_5\triangleq(\theta_1\cup\theta_2)\cap\theta_3\overset{\mathcal{M}_7}{\equiv}\emptyset$,
$\alpha_6\triangleq(\theta_1\cup\theta_3)\cap\theta_2\overset{\mathcal{M}_7}{\equiv}\emptyset$,
$\alpha_7\triangleq(\theta_2\cup\theta_3)\cap\theta_1\overset{\mathcal{M}_7}{\equiv}\emptyset $,
$ \alpha_8\triangleq\{(\theta_1\cap\theta_2)\cup\theta_3\} \cap(\theta_1\cup\theta_2)\overset{\mathcal{M}_7}{\equiv}\emptyset$ and
$\alpha_{11}\triangleq\theta_3\overset{\mathcal{M}_7}{\equiv}\emptyset$.
Therefore, one has now the following set of elements for $D^\Theta (\mathcal{M}_7(\Theta))$
\begin{equation*}
\begin{array}{|lll|}
\hline
& \text{Elements of } \; D^\Theta  \; \text{for} \; \mathcal{M}_7(\Theta)& \\
\hline
\alpha_0\triangleq \emptyset               & &                                 \\
\alpha_1\triangleq\theta_1\cap\theta_2\cap\theta_3 \overset{\mathcal{M}_7}{\equiv}\emptyset   & &    \\
\alpha_2\triangleq\theta_1\cap\theta_2 \overset{\mathcal{M}_7}{\equiv}\emptyset   & \alpha_3\triangleq\theta_1\cap\theta_3 \overset{\mathcal{M}_7}{\equiv}\emptyset  &    \alpha_4\triangleq\theta_2\cap\theta_3\overset{\mathcal{M}_7}{\equiv}\emptyset                          \\
\alpha_5\triangleq(\theta_1\cup\theta_2)\cap\theta_3\overset{\mathcal{M}_7}{\equiv}\emptyset   & \alpha_6\triangleq(\theta_1\cup\theta_3)\cap\theta_2\overset{\mathcal{M}_7}{\equiv}\emptyset &  \alpha_7\triangleq(\theta_2\cup\theta_3)\cap\theta_1\overset{\mathcal{M}_7}{\equiv}\emptyset  \\
  \alpha_8\triangleq\{(\theta_1\cap\theta_2)\cup\theta_3\} \cap(\theta_1\cup\theta_2)\overset{\mathcal{M}_7}{\equiv}\emptyset  & &  \\
\alpha_9\triangleq\theta_1 \neq \emptyset      & \alpha_{10}\triangleq\theta_2\neq\emptyset  & \alpha_{11}\triangleq\theta_3\overset{\mathcal{M}_7}{\equiv}\emptyset                                                   \\
\alpha_{12}\triangleq(\theta_1\cap\theta_2)\cup\theta_3\overset{\mathcal{M}_7}{\equiv}\emptyset     & \alpha_{13}\triangleq(\theta_1\cap\theta_3)\cup\theta_2\overset{\mathcal{M}_7}{\equiv} \alpha_{10}\neq   \emptyset &  \alpha_{14}\triangleq(\theta_2\cap\theta_3)\cup\theta_1\overset{\mathcal{M}_7}{\equiv}\alpha_{9}\neq\emptyset \\
\alpha_{15}\triangleq\theta_1\cup\theta_2\neq\emptyset    &      \alpha_{16}\triangleq\theta_1\cup\theta_3\overset{\mathcal{M}_7}{\equiv}\alpha_{9}\neq \emptyset  &                 \alpha_{17}\triangleq\theta_2\cup\theta_3\overset{\mathcal{M}_7}{\equiv}\alpha_{10}\neq \emptyset\\
\alpha_{18}\triangleq\theta_1\cup\theta_2\cup\theta_3\overset{\mathcal{M}_7}{\equiv}\alpha_{15}\neq \emptyset    & &  \\
\hline
\end{array}
\end{equation*}

$D^\Theta (\mathcal{M}_7(\Theta))$ reduces now to only four  different elements $\emptyset$, $\theta_1$, $\theta_2$, and $\theta_1\cup\theta_2$.

%************************************************
\section{DSm rule of combination for DSm hybrid models}
%*************************************************************

We present in this section a general DSm-hybrid rule of combination able to deal with any DSm hybrid models. We will show how this new general rule of combination works with all DSm hybrid models presented in the previous section and we list interesting properties of this new useful and powerful rule of combination.

\subsection{Notations}

Let $\Theta=\{\theta_1,\ldots\theta_n\}$ be a frame of {\it{partial discernment}} of the constrained fusion problem, and $D^\Theta$ the free distributive lattice (hyper-power set) generated by $\Theta$ and the empty set $\emptyset$ under $\cap$ and $\cup$ operators.   We need to distinguish between the empty set $\emptyset$, which belongs to $D^\Theta$, and by $\emptyset$ we understand a set which is empty all the time (we call it {\it{absolute emptiness}} or {\it{absolutely empty}}) independent of time, space and model, and all other sets from $D^\Theta$. For example $\theta_1\cap\theta_2$ or $\theta_1\cup\theta_2$ or only $\theta_i$ itself, $1\leq i \leq n$, etc, which could be or become empty at a certain time (if we consider a fusion dynamicity) or in a particular model $\mathcal{M}$ (but could not be empty in other model  and/or time) (we call a such element {\it{relative emptiness}} or {\it{relatively empty}}).  We’ll denote by $\boldsymbol{\emptyset}_{\mathcal{M}}$ the set of relatively empty such elements of $D^\Theta$ (i.e. which become empty in a particular model $\mathcal{M}$ or at a specific time).  $\boldsymbol{\emptyset}_{\mathcal{M}}$ is the set of integrity constraints which depends on the DSm model $\mathcal{M}$ under consideration, and the model $\mathcal{M}$ depends on the structure of its corresponding {\it{fuzzy Venn Diagram}} (number of elements in $\Theta$, number of non-empty intersections, and time in case of dynamic fusion). Through our convention $\emptyset\notin\boldsymbol{\emptyset}_{\mathcal{M}}$. 
Let’s note by $\boldsymbol{\emptyset}\triangleq \{\emptyset,\boldsymbol{\emptyset}_{\mathcal{M}}\}$ the set of all relatively and absolutely empty elements.\\

For any $A\in D^\Theta$, let $\phi(A)$ be the {\it{characteristic emptiness function}} of the set $A$, i.e. $\phi(A)= 1$ if  $A\notin \boldsymbol{\emptyset}$ and $\phi(A)= 0$ otherwise.  This function helps in assigning the value zero to all relatively or absolutely empty elements of $D^\Theta$ through the choice of DSm hybrid model $\mathcal{M}$.
Let's define the total ignorance on $\Theta = \{\theta_1,\theta_2, \ldots, \theta_n\}$ as
$I_t\triangleq\theta_1\cup \theta_2\cup \ldots\cup \theta_n$ and the set of relative ignorances as
$I_r\triangleq\{\theta_{i_1}\cup\ldots\cup\theta_{i_k}, \;\text{where}\; i_1, ..., i_k \in \{1,2,...,n\} \;  \text{and}\; 2\leq k\leq n-1\}$,
then the set of all kind of ignorances as $I = I_t \cup I_r$.
For any element $A$ in $D^\Theta$, one considers $u(A)$ as the union of all singletons $\theta_i$ that compose $A$.  For example, if $A$ is a singleton then $u(A)=A$; if  $A=\theta_1\cap \theta_2$ or $A=\theta_1\cup \theta_2$ then $u(A)=\theta_1\cup \theta_2$;
if $A=(\theta_1\cap \theta_2)\cup \theta_3$ then $u(A)=\theta_1\cup \theta_2\cup\theta_3$.  
; by convention $u(\emptyset)\triangleq\emptyset$.
The second summation of the DSm hybrid rule (see eq. \eqref{eq:DSmH2} and \eqref{eq:DSmHk} and denoted $S_2$ in the sequel) transfers the mass of $\emptyset$ [if any; sometimes, in rare cases, $m(\emptyset) > 0$ (for example in Ph. Smets' work); we want to catch this particular case as well] to the total ignorance $I_t=\theta_1\cup\theta_2\cup\ldots\cup\theta_n$.
The other part of the mass of relatively empty elements, $\theta_i$ and $\theta_j$ together for example, $i\neq j$, goes to the partial ignorance/uncertainty $m(\theta_i\cup\theta_j)$. $S_2$ multiplies, naturally following the DSm classic network architecture, only the elements of columns of absolutely and relatively empty sets, and then $S_2$ transfers the mass $m_1(X_1)m_2(X_2)\ldots m_k(X_k)$ either to the element $A\in D^\theta$ in the case when $A= u(X_1)\cup u(X_2)\cup \ldots \cup u(X_k)$ is not empty, or if $u(X1)\cup u(X2)\cup \ldots\cup u(X_k)$ is empty then the mass $m_1(X_1)m_2(X_2)… m_k(X_k)$ is transferred to the total ignorance.
We include all possible degenerate problems/models in this new DSmT hybrid framework, but the vacuous DSm-hybrid model $\mathcal{M}_{\emptyset}$ defined by the constraint $I_t= \theta_1\cup\theta_2\cup\ldots\cup\theta_n\overset{\mathcal{M}_{\emptyset}}{\equiv}\emptyset $ which is meaningless and useless.\\

We provide here the issue for programming the calculation of $u(X)$ from the binary representation of any proposition $X\in D^\Theta$ expressed in the Dezert-Smarandache's order \cite{Dezert_2003f,Dezert_Smarandache_2003}.
Let's consider the Smarandache's codification of elements $\theta_1,\ldots, \theta_n$.  One defines the anti-absorbing relationship as follows:
element $i$ anti-absorbs element $ij$ (with $i<j$), and let's use the notation $i << ij$, and also $j << ij$;
similarly $ij << ijk$ (with $i<j<k$), also $jk << ijk$ and $ik << ijk$.
This relationship is transitive, therefore $i << ij$ and $ij << ijk$ involve $i << ijk$; one can also write $i << ij << ijk$ as a chain;
similarly one gets $j << ijk$ and $k << ijk$.
The anti-absorbing relationship can be generalized for parts with any number of digits, i.e. when one uses the Smarandache  codification for the corresponding Venn diagram on $\Theta = \{\theta_1,\theta_2,\ldots, \theta_n\}$, with $n\geq 1$.
Between elements $ij$ and $ik$, or between $ij$ and $jk$ there is no anti-absorbing relationship, therefore the anti-absorbing relationship makes a partial order on the parts of the Venn diagram for the free DSm model.
If a proposition $X$ is formed by a part only, say $i_1 i_2\ldots i_r$, in the Smarandache codification, then $u(X)=\theta_{i_1}\cup \theta_{i_2}\cup\ldots\cup\theta_{i_r}$.
If $X$ is formed by two or more parts, the first step is to eliminate all anti-absorbed parts, ie. if $A << B$ then $u({A,B})=u(A)$;
generally speaking, a part $B$ is anti-absorbed by part $A$ if all digits of $A$ belong to $B$;  for an anti-absorbing chain $A_1 << A_2 << ... << A_s$ one takes $A_1$ only and the others are eliminated;
afterwards, when $X$ is anti-absorbingly irreducible, $u(X)$ will be the unions of all singletons whose indices occur in the remaining parts of $X$ - if one digit occurs many times it is taken only once.\\

See some examples for the case $n=3$:
$12 << 123$, i.e. $12$ anti-absorbs $123$.
Between $12$ and $23$ there is no anti-absorbing relationship.
\begin{itemize}
\item If $X=123$ then $u(X)=\theta_1\cup\theta_2\cup\theta_3$.
\item If $X=\{23,123\}$, then $23 << 123$, thus $u(\{23,123\})=u(23)$, because $123$ has been eliminated, hence $u(X)=u(23)=\theta_2\cup\theta_3$.
\item If $X=\{13,123\}$, then $13 << 123$, thus $u(\{13,123\})=u(13)=\theta_1\cup\theta_3$.
\item If $X=\{13,23,123\}$, then $13 << 123$, thus $u(\{13,23,123\})=u(\{13,23\})=\theta_1\cup\theta_2\cup\theta_3$ (one takes as theta indices each digit in the $\{13,23\}$) - if one digit is repeated it is taken only once; between $13$ and $23$ there is no relation of anti-absorbing.
\item If $X=\{3,13,23,123\}$, then $u(X)=u(\{3,13,23\})$ because $23 << 123$, then $u(\{3,13,23\})=u(\{3,13\})$ because $3 << 23$, then $u(\{3,13\})=u(3)=\theta_3$ because $3 << 13$.
\item If $X=\{1,12,13,23,123\}$, then one has the anti-absorbing chain: $1 << 12 << 123$, thus $u(X)=u(\{1,13,23\})=u(\{1,23\})$ because $1 << 13$, and finally $u(X)=\theta_1\cup\theta_2\cup\theta_3$.
\item If $X= \{1,2,12,13,23,123\}$, then $1<<12<<123$ and $2 << 23$ thus $u(X)=u(\{1,2,13\})=u(\{1,2\})$ because $1 << 13$, and finally $u(X)=\theta_1\cup \theta_2$.
\item If $X=\{2,12,3,13,23,123\}$, then $2 << 23 << 123$ and $3 << 13$ thus $u(X)=u(\{2,12,3\})$, but $2 << 12$ hence $u(X)=u(\{2,3\})=\theta_2\cup\theta_3$.
\end{itemize}

\clearpage
\newpage

\subsection{The DSm hybrid rule of combination for 2 sources}

To eliminate the {\it{degenerate vacuous fusion problem}} from the presentation, we assume from now on that the given DSm hybrid model $\mathcal{M}$ under consideration is always different from the vacuous model 
$\mathcal{M}_{\emptyset}$ (i.e. $I_t\neq\emptyset$). The DSm hybrid rule of combination, associated to a given DSm hybrid model $\mathcal{M}\neq\mathcal{M}_{\emptyset}$ , for two sources is defined for all $A\in D^\Theta$ as:

\begin{multline}
m_{\mathcal{M}(\Theta)}(A)\triangleq 
\phi(A)\Bigl[ \sum_{\overset{X_1,X_2\in D^\Theta}{(X_1\cap X_2)=A}} m_1(X_1)m_2(X_2)\\
+ \sum_{\overset{X_1,X_2\in\boldsymbol{\emptyset}}{[(u(X_1)\cup u(X_2))=A] \vee [(u(X_1)\cup u(X_2)\in\boldsymbol{\emptyset}) \wedge (A=I_t)]}}m_1(X_1)m_2(X_2)\\
+
\sum_{\overset{X_1,X_2\in D^\Theta}{\overset{(X_1\cup X_2)=A}{\overset{X_1\cap X_2\in\boldsymbol{\emptyset}}{}}}} m_1(X_1)m_2(X_2)
\Bigr]
 \label{eq:DSmH2}
\end{multline}

\noindent
The first sum entering in the previous formula corresponds to mass $m_{\mathcal{M}^f(\Theta)}(A)$ obtained by the classic DSm rule of combination \eqref{eq:DSMClassic2} based on the free-DSm model $\mathcal{M}^f$ (i.e. on the free lattice $D^\Theta$), i.e.
\begin{equation}
m_{\mathcal{M}^f(\Theta)}(A)\triangleq  \sum_{\overset{X_1,X_2\in D^\Theta}{X_1\cap X_2=A}} m_1(X_1)m_2(X_2)
\end{equation}

\noindent
The second sum entering in the formula of the DSm-hybrid rule of combination \eqref{eq:DSmH2} represents the mass of all relatively and absolutely empty sets which is transferred to the total or relative ignorances.\\

\noindent
The third sum entering in the formula of the DSm-hybrid rule of combination \eqref{eq:DSmH2} transfers the sum of relatively empty sets to the non-empty sets in the same way as it was calculated following the DSm classic rule.\\

\subsection{The DSm hybrid rule of combination for $k\geq 2$ sources}
%----------------------------------------------------------------------------------------------

The previous formula of DSm hybrid rule of combination can be generalized in the following way for all $A\in D^\Theta$ :
\begin{multline}
m_{\mathcal{M}(\Theta)}(A)\triangleq 
\phi(A)\Bigl[ \sum_{\overset{X_1,X_2,\ldots,X_k\in D^\Theta}{(X_1\cap X_2\cap\ldots\cap X_k)=A}} \prod_{i=1}^{k} m_i(X_i)\\
+ \sum_{\overset{X_1,X_2,\ldots,X_k\in\boldsymbol{\emptyset}}{[(u(X_1)\cup u(X_2)\cup \ldots \cup u(X_k))=A]\vee [(u(X_1)\cup u(X_2)\cup \ldots \cup u(X_k)\in\boldsymbol{\emptyset}) \wedge (A=I_t)]}} \prod_{i=1}^{k} m_i(X_i) \\
%
%+ \sum_{\overset{V\in\boldsymbol{\emptyset}, \: X_1,X_2,\ldots,X_k\in D^\Theta}{\overset{X_1\cup X_2\cup\ldots\cup X_k=A}{\overset{X_1\cap X_2\cap \ldots\cap X_k=V}{}}}}  
+ \sum_{\overset{X_1,X_2,\ldots,X_k\in D^\Theta}{\overset{(X_1\cup X_2\cup\ldots\cup X_k)=A}{\overset{X_1\cap X_2\cap \ldots\cap X_k\in\boldsymbol{\emptyset}}{}}}}  \prod_{i=1}^{k} m_i(X_i)
\Bigr]
 \label{eq:DSmHk}
\end{multline}

\noindent
The first sum entering in the previous formula corresponds to mass $m_{\mathcal{M}^f(\Theta)}(A)$ obtained by the classic DSm rule of combination \eqref{eq:DSMClassick} for $k$ sources of information based on the free-DSm model $\mathcal{M}^f$ (i.e. on the free lattice $D^\Theta$), i.e.
\begin{equation}
m_{\mathcal{M}^f(\Theta)}(A)\triangleq \sum_{\overset{X_1,X_2,\ldots,X_k\in D^\Theta}{(X_1\cap X_2\cap\ldots\cap X_k)=A}} \prod_{i=1}^{k} m_i(X_i)
\label{mfk}
\end{equation}

\subsection{Remark on the  DSm hybrid rule of combination}
%--------------------------------------------------------------------------

From \eqref{eq:DSmHk} and \eqref{mfk}, the previous general formula can be rewritten as 

\vspace{0.5cm}\begin{equation}
m_{\mathcal{M}(\Theta)}(A)\triangleq 
\phi(A)\Bigl[ S_1(A) + S_2(A) + S_3(A)\Bigr]
 \label{eq:DSmHkBis}
\end{equation}
\vspace{0.5cm}

\noindent
where
\begin{equation}
S_1(A)\equiv m_{\mathcal{M}^f(\Theta)}(A) \triangleq \sum_{\overset{X_1,X_2,\ldots,X_k\in D^\Theta}{(X_1\cap X_2\cap\ldots\cap X_k)=A}} \prod_{i=1}^{k} m_i(X_i)
\end{equation}
\vspace{1cm}
\begin{equation}
S_2(A)\triangleq \sum_{\overset{X_1,X_2,\ldots,X_k\in\boldsymbol{\emptyset}}{[(u(X_1)\cup u(X_2)\cup \ldots \cup u(X_k))=A]\vee [(u(X_1)\cup u(X_2)\cup \ldots \cup u(X_k)\in\boldsymbol{\emptyset}) \wedge (A=I_t)]}} \prod_{i=1}^{k} m_i(X_i)\end{equation}
\vspace{1cm}
\begin{equation}
S_3(A)\triangleq\sum_{\overset{X_1,X_2,\ldots,X_k\in D^\Theta}{\overset{(X_1\cup X_2\cup\ldots\cup X_k)=A}{\overset{X_1\cap X_2\cap \ldots\cap X_k\in\boldsymbol{\emptyset}}{}}}}  \prod_{i=1}^{k} m_i(X_i)
\end{equation}
\vspace{1cm}

\noindent
and thus, this combination can be viewed actually as a two steps procedure as follows:
\begin{itemize}
\item Step 1:  Evaluate the combination of the sources over the free lattice $D^\Theta$ by the classical DSm rule of combination to get for all $A\in D^\Theta$, $S_1(A)=m_{\mathcal{M}^f(\Theta)}(A)$ using \eqref{mfk}. This step preserves the commutativity and associativity properties of the combination.\\

\item Step 2 : Transfer the masses of the integrity constraints of the DSm hybrid model $\mathcal{M}$ according to formula \eqref{eq:DSmHkBis}. Note that this step is necessary only if one has reliable information about the real constraints involved in the fusion problem under consideration.
\end{itemize}

The second step does not preserve the associativity property but this is not a fundamental requirement in most of fusion systems actually.
If one really wants to preserve optimality of the fusion rule, one has first to combine all sources using classical DSm rule (with any clustering of sources) and the ultimate step will consist to adapt basic belief masses according to the integrity constraints of the model $\mathcal{M}$.\\

 If one first adapts the local basic belief masses $m_1(.),...m_k(.)$ to the hybrid-model $\mathcal{M}$ and afterwards one applies the combination rules, the fusion result becomes only suboptimal because some information is lost during the transfer of masses of integrity constraints. The same remark holds if the transfer of masses of integrity constraints is done at some intermediate steps after the fusion of $m$ sources with $m<k$.\\
 
Let's note also that this formula of transfer is more general (because we include the possibilities to introduce both exclusivity constraints and non-existential constraints as well) and more precise (because we explicitly consider all different relative emptiness of elements into the general transfer formula \eqref{eq:DSmHkBis}) than the generic transfer formulas used in the DST framework proposed as alternative rules  to the Dempster's rule of combination \cite{Lefevre_2002} and discussed in section 5.9.

\clearpage
\newpage

%--------------------------------------------------------
\subsection{Property of the DSm Hybrid Rule}
%--------------------------------------------------------

\begin{equation}
\sum_{A\in D^\Theta} m_{\mathcal{M}(\Theta)}(A) = \sum_{A\in D^\Theta} \phi(A)\Bigl[ S_1(A) + S_2(A) + S_3(A)\Bigr]=1
\end{equation}
\vspace{0.5cm}

{\bf{Proof}}: Let's first prove that $\sum_{A\in D^\Theta} m(A)=1$ where all masses $m(A)$ are obtained by the DSm classic rule. Let's consider each mass $m_i(.)$ provided by the $i$th source of information, for $1\leq i\leq k$, as a vector of $d={\mid D^\Theta\mid}$ dimension, whose sum of components is equal to one, i.e. $m_i(D^\Theta)=[m_{i1},m_{i2},\ldots,m_{id}]$, and $\sum_{j=1,d} m_{ij}=1$. Thus, for $k\geq 2$ sources of information, the mass matrix becomes
\vspace{0.5cm}
\begin{equation*}
\mathbf{M}=
\begin{bmatrix}
m_{11} & m_{12} & \ldots & m_{1d}\\
m_{21} & m_{22} & \ldots & m_{2d}\\
 \ldots &  \ldots & \ldots & \ldots\\
m_{k1} & m_{k2} & \ldots & m_{kd}\\
\end{bmatrix}
\end{equation*}

\vspace{0.5cm}

\noindent If one notes the sets in $D^\Theta$ by $A_1$, $A_2$, ..., $A_d$ (doesn't matter in what order one lists them) then the column $(j)$ in the matrix represents the masses assigned to $A_j$ by each source of information $s_1$, $s_2$, $\ldots$, $s_k$; for example $s_i(A_j)=m_{ij}$, where $1\leq i\leq k$.
According to the DSm network architecture \cite{Dezert_2003f}, all the products in this network will have the form $m_{1j_1}m_{2j_2}\ldots m_{kj_k}$, i.e. one element only from each matrix row, and no restriction about the number of elements from each matrix column, $1\leq j_1,j_2,\ldots,j_k\leq d$. Each such product will enter in the fusion mass of one set only from $D^\Theta$. Hence the sum of all components of the fusion mass is equal to the sum of all these products, which is equal to 
\vspace{0.5cm}
\begin{equation}
\prod_{i=1}^{k} \sum_{j=1}^{d} m_{ij}=\prod_{i=1}^{k} 1 = 1
\label{eq:Formula12}
\end{equation}
\vspace{0.5cm}
The DSm hybrid rule has three sums $S_1$, $S_2$, and $S_3$. Let's separate the mass matrix $\mathbf{M}$ into two disjoint sub-matrices $\mathbf{M}_{\emptyset}$ formed by the columns of all absolutely and relatively empty sets, and $\mathbf{M}_N$ formed by the columns of all non-empty sets. According to the DSm network architecture (for $k\geq 2$ rows):
\begin{itemize}
\item $S_1$ is the sum of all products resulted from the multiplications of the columns of $\mathbf{M}_N$ following the DSm network architecture such that the intersection of their corresponding sets is non-empty, i.e. the sum of masses of all non-empty sets before any mass of absolutely or relatively empty sets could be transferred to them;
\item $S_2$ is the sum of all products resulted from the multiplications of $\mathbf{M}_{\emptyset}$ following the DSm network architecture, i.e. a partial sum of masses of absolutely and relatively empty sets transferred to the ignorances in $I\triangleq I_t\cup I_r$ or to singletons of $\Theta$.
\item 
$S_3$ is the sum of all the products resulted from the multiplications of the columns of $\mathbf{M}_N$ and $\mathbf{M}_{\emptyset}$ together, following the DSm network architecture, but such that at least a column is from each of them, and also the sum of all products of columns of $\mathbf{M}_N$  such that the intersection of their corresponding sets is empty (what did not enter into the previous sum $S_1$), i.e. the remaining sum of masses of absolutely or relatively empty sets transferred to the non-empty sets of the DSm hybrid model $\mathcal{M}$.
\end{itemize}
If one now considers all the terms (each such term is a product of the form   $m_{1j_1}m_{2j_2}\ldots m_{kj_k}$) of these three sums, we get exactly the same terms as in the DSm network architecture for the DSm classic rule, thus the sum of all terms occurring in $S_1$, $S_2$, and $S_3$ is 1 (see formula \eqref{eq:Formula12}) which completes the proof. DSm hybrid rule naturally derives from the DSm classic rule.
\clearpage
\newpage
Entire masses of relatively and absolutely empty sets in a given DSm hybrid model $\mathcal{M}$ are transferred to non-empty sets according to the above and below formula  \eqref{eq:DSmHkBis} and thus
\begin{equation}
\forall A\in \boldsymbol{\emptyset} \subset D^\Theta, \quad m_{\mathcal{M}(\Theta)}(A)=0
\end{equation}
The entire mass of a  relatively empty set (from $D^\Theta$) which has in its expression $\theta_{j_1}$, $\theta_{j_2}$, $\ldots$, $\theta_{j_r}$, with $1\leq r\leq n$ will generally be distributed among the $\theta_{j_1}$, $\theta_{j_2}$, $\ldots$, $\theta_{j_r}$ or their unions or intersections, and the distribution follows the way of multiplication from the DSm classic rule, explained by the DSm network architecture \cite{Dezert_2003f}. Thus, because nothing is lost, nothing is gained, the sum of all $m_{\mathcal{M}(\Theta)}(A)$ is equal to 1 as just proved previously, and fortunately no normalization constant is needed which could bring a lost of information in the fusion rule.\\

The three summations $S_1(.)$, $S_3(.)$ and $S_3(.)$ are disjoint because:
\begin{itemize}
\item $S_1(.)$ multiplies the columns corresponding to non-emptysets only - but such that the intersections of the sets corresponding to these columns are non-empty [following the definition of DSm classic rule];
\item $S_2(.)$ multiplies the columns corresponding to absolutely and relatively emptysets only;
\item $S_3(.)$ multiplies: 
\begin{itemize}
\item[a)] either the columns corresponding to absolutely or relatively emptysets with the columns corresponding to non-emptysets such that at least a column corresponds to an absolutely or relatively emptyset and at least a column corresponds to a non-emptyset, 
\item[b)]  or the columns corresponding to non-emptysets - but such that the intersections of the sets corresponding to these columns are empty.
\end{itemize}
\end{itemize}
The multiplications are following the DSm network architecture, i.e. any product has the above general form: $m_{1j_1}m_{2j_2}\ldots m_{kj_k}$, i.e.
any product contains as factor one element only from each row of the mass matrix $\mathbf{M}$ and the total number of factors in a product is equal to $k$.
The function $\phi(A)$ automatically assigns the value zero to the mass of any empty set, and allows the calculation of masses of all non-emptysets.\\

%------------------------------------------------------------------------------
\subsection{On the programming of the DSm hybrid rule}
%------------------------------------------------------------------------------

We briefly give here  an issue for a fast programming of DSm rule of combination. Let's consider $\Theta = \{\theta_1, \theta_2, \ldots, \theta_n\}$, the sources $\mathcal{B}_1$, $\mathcal{B}_2$,$\ldots$, $\mathcal{B}_k$, and $p= \min\{n, k\}$. One needs to check only the focal sets, i.e. sets (i.e. propositions) whose masses assigned to them by these sources are not all zero. Thus, if $\mathbf{M}$ is the mass matrix, and we consider a set $A_j$ in $D^\Theta$, then the column $(j)$ corresponding to $A_j$, i.e. $(m_{1j}  \; m_{2j}\; \ldots \; m_{kj})$ transposed has not to be identical to the null-vector of $k$-dimension $(0\;  0 \; \ldots \; 0)$ transposed.
Let $D^\Theta(\text{step}_1)$ be formed by all focal sets at the beginning (after sources $\mathcal{B}_1$, $\mathcal{B}_2$,$\ldots$, $\mathcal{B}_k$ have assigned massed to the sets in $D^\Theta$). Applying the DSm classic rule, besides the sets in $D^\Theta(\text{step}_1)$ one adds $r$-intersections of sets in $D^\Theta(\text{step}_1)$, thus:

$$D^\Theta(\text{step}_2) = D^\Theta(\text{step}_1) \vee \{A_{i_1}\wedge A_{i_2}\wedge\ldots\wedge A_{i_r}\}$$

\noindent where $A_{i_1}$, $A_{i_2}$, \ldots, $A_{i_r}$ belong to $D^\Theta(\text{step}_1)$ and $2\leq r\leq p$.\\

Applying the DSm hybrid rule, due to its $S_2$ and $S_3$ summations, besides the sets in $D^\Theta(\text{step}_2)$ one adds $r$-unions of sets and the total ignorance in $D^\Theta(\text{step}_2)$, thus:

$$D^\Theta(\text{step}_3) = D^\Theta(\text{step}_2) \vee I_t \vee \{A_{i_1}\vee A_{i_2}\vee\ldots\vee A_{i_r}\}$$

\noindent where $A_{i_1}$, $A_{i_2}$, \ldots, $A_{i_r}$ belong to $D^\Theta(\text{step}_2)$ and $2\leq r\leq p$.\\

This means that instead of computing the masses of all sets in $D^\Theta$, one needs to first compute the masses of all focal sets (step 1), second the masses of their $r$-intersections (step 2), and third the masses of $r$-unions of all previous sets and the mass of total ignorance (step 3).

\clearpage
\newpage
%--------------------------------------------------------------------------------------------
\subsection{Application of the DSm Hybrid rule on previous examples}
%--------------------------------------------------------------------------------------------

We present in this section some numerical results of the DSm hybrid rule of combination for 2 independent sources of information. We examine the seven previous examples in order to help the reader to check by himself or herself the validity of our new general formula. Due to space limitation, we will not go in details on all the derivations steps, we will just present main intermediary results (i.e. the value of the three summations) involved into the general formula \eqref{eq:DSmH2}. The results have been first obtained by hands and then be validated by MatLab programming. We denote
$$S_1(A)\equiv m_{\mathcal{M}^f(\Theta)}(A)\triangleq  \sum_{\overset{X_1,X_2\in D^\Theta}{X_1\cap X_2=A}} m_1(X_1)m_2(X_2)$$
$$S_2(A)\triangleq \sum_{\overset{X_1,X_2\in\boldsymbol{\emptyset}}{[u(X_1)\cup u(X_2)=A] \vee [(u(X_1)\cup u(X_2)\in\boldsymbol{\emptyset}) \wedge (A=I_t)]}}m_1(X_1)m_2(X_2)$$
$$S_3(A)\triangleq \sum_{\overset{X_1,X_2\in D^\Theta}{\overset{X_1\cup X_2=A}{\overset{X_1\cap X_2\in\boldsymbol{\emptyset}}{}}}} m_1(X_1)m_2(X_2)$$
Now let consider $\Theta=\{\theta_1,\theta_2,\theta_3\}$ and the two following independent bodies of evidence $\mathcal{B}_1$ and $\mathcal{B}_2$ with the generalized basic belief assignments\footnote{A general numerical example with $m_1(A)>0$ and $m_2(A)>0$ for all $A\neq\emptyset\in D^\Theta$ will be briefly presented in next section.} $m_1(.)$ and $m_2(.)$ given in the following table\footnote{The order of elements of $D^\Theta$ corresponds here to the order obtained from the generation of isotone Boolean functions - see \cite{Dezert_2003f} for details.}. The right column of the table indicates the result of the fusion obtained by the classical DSm rule of combination.
\begin{equation*}
\begin{array}{llll}
\hline
\text{Element} \: A \: \text{of} \; D^\Theta&   m_1(A) & m_2(A)  & m_{\mathcal{M}^f(\Theta)}(A) \\
\hline
\emptyset               &    0 &     0               &  0             \\
\theta_1\cap\theta_2\cap\theta_3 & 0 &   0 & 0.16 \\
\theta_2\cap\theta_3 & 0 &0.20        & 0.19                \\
\theta_1\cap\theta_3 &  0.10 & 0 & 0.12\\
(\theta_1\cup\theta_2)\cap\theta_3  & 0 & 0 & 0.01\\
\theta_3 & 0.30 & 0.10 & 0.10\\                          
\theta_1\cap\theta_2  & 0.10 & 0.20 & 0.22\\
(\theta_1\cup\theta_3)\cap\theta_2 &  0 & 0 & 0.05\\
(\theta_2\cup\theta_3)\cap\theta_1& 0 & 0 & 0 \\
\{(\theta_1\cap\theta_2)\cup\theta_3\} \cap(\theta_1\cup\theta_2) & 0 & 0& 0 \\
(\theta_1\cap\theta_2)\cup\theta_3 & 0 & 0 & 0\\
\theta_2 & 0.20 & 0.10 & 0.03\\                          
(\theta_1\cap\theta_3)\cup\theta_2 &  0 & 0 & 0\\
\theta_2\cup\theta_3 & 0 & 0 & 0\\
\theta_1 & 0.10 & 0.20 & 0.08\\
(\theta_2\cap\theta_3)\cup\theta_1& 0 & 0 & 0.02\\
\theta_1\cup\theta_3 &  0.10 & 0.20 & 0.02\\
\theta_1\cup\theta_2 & 0.10 & 0 & 0\\
\theta_1\cup\theta_2\cup\theta_3 & 0 & 0 & 0
\end{array}
\end{equation*}
The following subsections present the numerical results obtained by the DSm hybrid rule on the seven previous examples.
The tables show all the values of $\phi(A)$, $S_1(A)$, $S_2(A)$ and $S_3(A)$ to help the reader to check by himself or herself the validity of these results.
It is important to note that the values of $S_1(A)$, $S_2(A)$ and $S_3(A)$ when $\phi(A)=0$ do not need to be  computed in practice but are provided here only for a checking purpose.

\clearpage
\newpage

\subsubsection{Application of the DSm Hybrid rule on example 1}
%---------------------------------------------------------------------------------

Here is the numerical result corresponding to example 1 with the hybrid-model $\mathcal{M}_1$ (i.e with the exclusivity constraint $\theta_1\cap\theta_2\cap\theta_3  \overset{\mathcal{M}_1}{\equiv}\emptyset$). The right column of the table provides the result obtained using the DSm hybrid rule,
ie. $\forall A\in D^\Theta$, $$m_{\mathcal{M}_1(\Theta)}(A)=\phi(A)\bigl[ S_1(A) +  S_2(A) + S_3(A)\bigr] $$

\begin{equation*}
\begin{array}{|l|rrrr|r|}
\hline
\text{Element} \: A \: \text{of} \; D^\Theta&   \phi(A) & S_1(A) &  S_2(A)& S_3(A) & m_{\mathcal{M}_1(\Theta)}(A)\\
\hline
\emptyset               & 0 & 0 & 0 & 0 & 0\\
\theta_1\cap\theta_2\cap\theta_3  \overset{\mathcal{M}_1}{\equiv}\emptyset & 0 & 0.16 & 0 & 0 & 0\\
\theta_2\cap\theta_3 & 1 & 0.19 & 0 & 0 & 0.19\\
\theta_1\cap\theta_3 & 1 & 0.12 & 0 & 0 & 0.12\\ 
(\theta_1\cup\theta_2)\cap\theta_3  & 1 & 0.01  & 0 & 0 .02 & 0.03\\ 
\theta_3 & 1 & 0.10  & 0 & 0  & 0.10\\      
\theta_1\cap\theta_2  & 1 & 0.22  & 0 & 0  & 0.22\\
(\theta_1\cup\theta_3)\cap\theta_2 & 1 & 0.05  & 0 & 0.02  & 0.07\\
(\theta_2\cup\theta_3)\cap\theta_1& 1 & 0  & 0 & 0.02  & 0.02\\
\{(\theta_1\cap\theta_2)\cup\theta_3\} \cap(\theta_1\cup\theta_2) & 1 & 0  & 0 & 0  & 0\\
(\theta_1\cap\theta_2)\cup\theta_3 & 1 & 0  & 0 & 0.07  & 0.07\\  
\theta_2 & 1 & 0.03 & 0 & 0  & 0.03\\             
(\theta_1\cap\theta_3)\cup\theta_2 & 1 & 0 & 0 & 0.01  & 0.01\\
\theta_2\cup\theta_3 & 1 & 0 & 0 & 0  & 0\\
\theta_1 & 1 & 0.08 & 0 & 0  & 0.08\\
(\theta_2\cap\theta_3)\cup\theta_1& 1 & 0.02 & 0 & 0.02  & 0.04\\
\theta_1\cup\theta_3 & 1 & 0.02 & 0 & 0  & 0.02\\
\theta_1\cup\theta_2 & 1 & 0 & 0 & 0  & 0\\
\theta_1\cup\theta_2\cup\theta_3 & 1 & 0 & 0 & 0  & 0\\
\hline
\end{array}
\quad
\mathbf{D}_{\mathcal{M}_1}=\begin{bmatrix}
    0  &   0    & 0   &  0   &  0   &  0   \\
    0  &   0   &  0   &  0  &   0  &   1   \\
     0  &   0  &   0  &   0  &   1  &   0   \\
     0  &   0  &   0  &   0  &   1  &   1  \\
     0   &  0  &   0  &   1  &   1  &   1   \\
     0  &   0   &  1   &  0  &   0  &   0   \\
     0  &   0  &   1  &   0  &   0  &   1   \\
     0  &   0  &   1  &   0  &   1  &   0  \\
     0  &   0   &  1  &   0  &   1  &   1  \\
     0  &   0   &  1  &   1   &  1  &   1  \\
     0  &   1  &   1 &    0  &   0   &  1  \\
     0  &   1  &   1  &   0  &   1   &  1  \\
     0  &   1  &   1  &   1   &  1  &   1  \\
     1  &   0  &   1  &   0  &   1  &   0 \\
     1  &   0  &   1  &   0  &   1  &   1  \\
     1  &   0  &   1  &   1  &   1  &   1 \\
     1   &  1  &   1  &   0   &  1  &   1  \\
     1  &   1  &   1  &   1  &   1  &   1  
\end{bmatrix}
\end{equation*}

From the previous table of this first numerical example, we see in column corresponding to $S_3(A)$ how the initial combined mass $m_{\mathcal{M}^f(\Theta)}(\theta_1\cap\theta_2\cap\theta_3)\equiv S_1(\theta_1\cap\theta_2\cap\theta_3)=0.16$ is transferred (due to the constraint of $\mathcal{M}_1$) only onto the elements $(\theta_1\cup\theta_2)\cap\theta_3$, $(\theta_1\cup\theta_3)\cap\theta_2$, $(\theta_2\cup\theta_3)\cap\theta_1$, $(\theta_1\cap\theta_2)\cup\theta_3$, $(\theta_1\cap\theta_3)\cup\theta_2$, and $(\theta_2\cap\theta_3)\cup\theta_1$  of $D^\Theta$. We can easily check that the sum of the elements of the column for $S_3(A)$ is equal to $m_{\mathcal{M}^f(\Theta)}(\theta_1\cap\theta_2\cap\theta_3)=0.16$ as required.
Thus after introducing the constraint, the initial hyper-power set $D^\Theta$ reduces to 18 elements as follows
\begin{equation*}
\begin{split}
D^\Theta_{\mathcal{M}_1}=\{\emptyset,\theta_2\cap\theta_3,\theta_1\cap\theta_3,(\theta_1\cup\theta_2)\cap\theta_3,\theta_3,\theta_1\cap\theta_2,(\theta_1\cup\theta_3)\cap\theta_2,(\theta_2\cup\theta_3)\cap\theta_1,\{(\theta_1\cap\theta_2)\cup\theta_3\} \cap(\theta_1\cup\theta_2),\\
(\theta_1\cap\theta_2)\cup\theta_3,\theta_2,(\theta_1\cap\theta_3)\cup\theta_2,\theta_2\cup\theta_3,\theta_1,(\theta_2\cap\theta_3)\cup\theta_1,\theta_1\cup\theta_3,\theta_1\cup\theta_2,\theta_1\cup\theta_2\cup\theta_3\}
\end{split}
\end{equation*}

\noindent
As detailed in \cite{Dezert_2003f}, the elements of $D^\Theta_{\mathcal{M}_1}$ can be described and encoded by the matrix product $\mathbf{D}_{\mathcal{M}_1}\cdot\mathbf{u}_{\mathcal{M}_1}$ with $\mathbf{D}_{\mathcal{M}_1}$ given above
and the basis vector $\mathbf{u}_{\mathcal{M}_1}$ defined as
$\mathbf{u}_{\mathcal{M}_1}=[<1>    <2>    <12>     <3>    <13>    <23>]'$. Actually $\mathbf{u}_{\mathcal{M}_1}$ is directly obtained from $\mathbf{u}_{\mathcal{M}^f}$\footnote{$\mathbf{D}_{\mathcal{M}^f}$ was denoted $\mathbf{D}_n$ and $\mathbf{u}_{\mathcal{M}^f}$ as $\mathbf{u}_n$ in reference \cite{Dezert_2003f}.} by removing its component $<123>$ corresponding to the constraint introduced by the model $\mathcal{M}_1$.\\

In general, the encoding matrix $\mathbf{D}_{\mathcal{M}}$ for a given DSm hybrid model $\mathcal{M}$ is obtained from $\mathbf{D}_{\mathcal{M}^f}$ by removing all its columns corresponding to the constraints of the chosen model $\mathcal{M}$ and all the rows corresponding to redundant/equivalent propositions. In this particular example with model $\mathcal{M}_1$, we will just have to remove the last column of $\mathbf{D}_{\mathcal{M}^f}$ to get  $\mathbf{D}_{\mathcal{M}_1}$ and no row is removed from $\mathbf{D}_{\mathcal{M}^f}$ because there is no redundant/equivalent proposition involved in this example. This suppression of some rows of $\mathbf{D}_{\mathcal{M}^f}$ will however occur in next examples. We encourage the reader to consult the references \cite{Dezert_2003f,Dezert_Smarandache_2003} for explanations and details about the generation, the encoding and the partial ordering of hyper-power sets.

\clearpage
\newpage

%*********************************************************************
\subsubsection{Application of the DSm Hybrid rule on example 2}
%*********************************************************************
 
Here is the numerical result corresponding to example 2 with the hybrid-model $\mathcal{M}_2$ (i.e with the exclusivity constraint $\theta_1\cap\theta_2\overset{\mathcal{M}_2}{\equiv}\emptyset \Rightarrow \theta_1\cap\theta_2\cap\theta_3\overset{\mathcal{M}_2}{\equiv}\emptyset$). One gets now

\begin{equation*}
\begin{array}{|l|rrrr|r|}
\hline
\text{Element} \: A \: \text{of} \; D^\Theta&   \phi(A) & S_1(A) & S_2(A)& S_3(A) & m_{\mathcal{M}_2(\Theta)}(A)\\
\hline
\emptyset               & 0 & 0  & 0 & 0 & 0\\
\theta_1\cap\theta_2\cap\theta_3  \overset{\mathcal{M}_2}{\equiv}\emptyset & 0 & 0.16  & 0 & 0 & 0\\
\theta_2\cap\theta_3 & 1 & 0.19  & 0 & 0 & 0.19\\
\theta_1\cap\theta_3 & 1 & 0.12  & 0 & 0 & 0.12\\ 
(\theta_1\cup\theta_2)\cap\theta_3  & 1 & 0.01 & 0 & 0 .02 & 0.03\\ 
\theta_3 & 1 & 0.10 & 0 & 0  & 0.10\\      
% ii=6 ==> atom Aii=0  0  1  0  0  0  1
\theta_1\cap\theta_2\overset{\mathcal{M}_2}{\equiv}\emptyset   & 0 & 0.22 & 0 & 0.02  & 0\\
(\theta_1\cup\theta_3)\cap\theta_2 \overset{\mathcal{M}_2}{\equiv}\theta_2\cap\theta_3 & 1 & 0.05 & 0 & 0.02  & 0.07\\
(\theta_2\cup\theta_3)\cap\theta_1\overset{\mathcal{M}_2}{\equiv}\theta_1\cap\theta_3 & 1 & 0 & 0 & 0.02  & 0.02\\
\{(\theta_1\cap\theta_2)\cup\theta_3\} \cap(\theta_1\cup\theta_2) \overset{\mathcal{M}_2}{\equiv}(\theta_1\cup\theta_2)\cap\theta_3
& 1 & 0 & 0 & 0  & 0\\
(\theta_1\cap\theta_2)\cup\theta_3 \overset{\mathcal{M}_2}{\equiv} \theta_3 & 1 & 0 & 0 & 0.07  & 0.07\\  
\theta_2 & 1 & 0.03 & 0 & 0.05  & 0.08\\  
(\theta_1\cap\theta_3)\cup\theta_2 & 1 & 0 & 0 & 0.01  & 0.01\\
\theta_2\cup\theta_3 & 1 & 0 & 0 & 0  & 0\\
\theta_1 & 1 & 0.08 & 0 & 0.04  & 0.12\\
(\theta_2\cap\theta_3)\cup\theta_1& 1 & 0.02 & 0 & 0.02  & 0.04\\
\theta_1\cup\theta_3 & 1 & 0.02 & 0 & 0.04 & 0.06\\
\theta_1\cup\theta_2 & 1 & 0 & 0.02 & 0.07  & 0.09\\
\theta_1\cup\theta_2\cup\theta_3 & 1 & 0 & 0 & 0  & 0\\
\hline
\end{array}
\end{equation*}

From the previous table of this numerical example, we see in column corresponding to $S_3(A)$ how the initial combined masses $m_{\mathcal{M}^f(\Theta)} (\theta_1\cap\theta_2\cap\theta_3) \equiv S_1(\theta_1\cap\theta_2\cap\theta_3)=0.16$ and 
$m_{\mathcal{M}^f(\Theta)}(\theta_1\cap\theta_2)\equiv S_1(\theta_1\cap\theta_2)=0.22$
are transferred (due to the constraint of $\mathcal{M}_2$) onto some elements of $D^\Theta$. 
We can easily check that the sum of the elements of the column for $S_3(A)$ is equal to $0.16+0.22=0.38$. Because some elements of $D^\Theta$ are now equivalent due to the constraints of $\mathcal{M}_2$, we have to sum all the masses corresponding to same equivalent propositions/elements (by example $\{(\theta_1\cap\theta_2)\cup\theta_3\} \cap(\theta_1\cup\theta_2) \overset{\mathcal{M}_2}{\equiv}(\theta_1\cup\theta_2)\cap\theta_3$). This can be viewed as the final compression step. One then gets the reduced hyper-power set $D^\Theta_{\mathcal{M}_2}$ having now 13 different elements with the combined belief masses presented in the following table. The basis vector $\mathbf{u}_{\mathcal{M}_2}$ and the encoding matrix $\mathbf{D}_{\mathcal{M}_2}$ for the elements of $D^\Theta_{\mathcal{M}_2}$ are given by $\mathbf{u}_{\mathcal{M}_2}=[<1>    <2>    <3>    <13>    <23>]'$ and below. Actually $\mathbf{u}_{\mathcal{M}_2}$ is directly obtained from $\mathbf{u}_{\mathcal{M}^f}$ by removing its components $<12>$ and $<123>$ corresponding to the constraints introduced by the model $\mathcal{M}_2$.
\begin{equation*}
\begin{array}{|l|r|}
\hline
\text{Element} \: A \: \text{of} \; D^\Theta_{\mathcal{M}_2} & m_{\mathcal{M}_2(\Theta)}(A)\\
\hline
\emptyset               & 0\\
\theta_2\cap\theta_3 &  0.19+0.07=0.26\\
\theta_1\cap\theta_3 & 0.12+0.02=0.14\\ 
(\theta_1\cup\theta_2)\cap\theta_3  & 0.03+0=0.03\\ 
\theta_3 & 0.10+0.07=0.17\\      
\theta_2 & 0.08\\  
(\theta_1\cap\theta_3)\cup\theta_2 & 0.01\\
\theta_2\cup\theta_3 & 0\\
\theta_1 & 0.12\\
(\theta_2\cap\theta_3)\cup\theta_1& 0.04\\
\theta_1\cup\theta_3 & 0.06\\
\theta_1\cup\theta_2 & 0.09\\
\theta_1\cup\theta_2\cup\theta_3 &  0\\
\hline
\end{array}
\qquad\text{and}\qquad
\mathbf{D}_{\mathcal{M}_2}=\begin{bmatrix}
    0  &   0       &  0   &  0   &  0   \\    
    0  &   0      &  0  &   0  &   1   \\   
     0  &   0    &   0  &   1  &   0   \\   
     0  &   0    &   0  &   1  &   1  \\
     0   &  0    &   1  &   1  &   1   \\
     0  &   1   &    0  &   0   &  1  \\
     0  &   1    &   0  &   1   &  1  \\
     0  &   1    &   1   &  1  &   1  \\
     1  &   0    &   0  &   1  &   0 \\
     1  &   0    &   0  &   1  &   1  \\
     1  &   0    &   1  &   1  &   1 \\
     1   &  1    &   0   &  1  &   1  \\
     1  &   1    &   1  &   1  &   1  
\end{bmatrix}
\end{equation*}

\clearpage
\newpage

\subsubsection{Application of the DSm Hybrid rule on example 3}
%---------------------------------------------------------------------------------

Here is the numerical result corresponding to example 3 with the hybrid-model $\mathcal{M}_3$ (i.e with the exclusivity constraint $(\theta_1\cup\theta_3)\cap\theta_2 \overset{\mathcal{M}_3}{\equiv}\emptyset$). This constraint implies directly $\theta_1\cap\theta_2\cap\theta_3\overset{\mathcal{M}_3}{\equiv}\emptyset$, $\theta_1\cap\theta_2\overset{\mathcal{M}_3}{\equiv}\emptyset$ and $\theta_2\cap\theta_3\overset{\mathcal{M}_3}{\equiv}\emptyset$. One gets now
\begin{equation*}
\begin{array}{|l|rrrrr|r|}
\hline
\text{Element} \: A \: \text{of} \; D^\Theta&   \phi(A) & S_1(A) & S_2(A)& S_3(A) & m_{\mathcal{M}_3(\Theta)}(A)\\
\hline
\emptyset               & 0 & 0 & 0 & 0 & 0\\
\theta_1\cap\theta_2\cap\theta_3  \overset{\mathcal{M}_3}{\equiv}\emptyset & 0 & 0.16 & 0 & 0 & 0\\
\theta_2\cap\theta_3\overset{\mathcal{M}_3}{\equiv}\emptyset & 0 & 0.19 & 0 & 0 & 0\\
\theta_1\cap\theta_3 & 1 & 0.12 & 0 & 0 & 0.12\\ 
(\theta_1\cup\theta_2)\cap\theta_3  \overset{\mathcal{M}_3}{\equiv} \theta_1\cap\theta_3 & 1 & 0.01 & 0 & 0 .02 & 0.03\\ 
\theta_3 & 1 & 0.10 & 0 & 0.06  & 0.16\\      
\theta_1\cap\theta_2\overset{\mathcal{M}_3}{\equiv}\emptyset   & 0 & 0.22 & 0 & 0.02  & 0\\
(\theta_1\cup\theta_3)\cap\theta_2 \overset{\mathcal{M}_3}{\equiv}\emptyset & 0 & 0.05 & 0 & 0.02  & 0\\
(\theta_2\cup\theta_3)\cap\theta_1\overset{\mathcal{M}_3}{\equiv} \theta_1\cap\theta_3 & 1 & 0 & 0 & 0.02  & 0.02\\
\{(\theta_1\cap\theta_2)\cup\theta_3\} \cap(\theta_1\cup\theta_2) \overset{\mathcal{M}_3}{\equiv} \theta_1\cap\theta_3
& 1 & 0 & 0 & 0  & 0\\
(\theta_1\cap\theta_2)\cup\theta_3 \overset{\mathcal{M}_3}{\equiv}  \theta_3 & 1 & 0 & 0 & 0.07  & 0.07\\  
\theta_2 & 1 & 0.03 & 0 & 0.09  & 0.12\\  
(\theta_1\cap\theta_3)\cup\theta_2 & 1 & 0 & 0 & 0.01  & 0.01\\
\theta_2\cup\theta_3 & 1 & 0 & 0 & 0.05  & 0.05\\
\theta_1 & 1 & 0.08 & 0 & 0.04  & 0.12\\
(\theta_2\cap\theta_3)\cup\theta_1 \overset{\mathcal{M}_3}{\equiv}  \theta_1 & 1 & 0.02 & 0 & 0.02  & 0.04\\
\theta_1\cup\theta_3 & 1 & 0.02 & 0 & 0.06 & 0.08\\
\theta_1\cup\theta_2 & 1 & 0 & 0.02 & 0.09  & 0.11\\
\theta_1\cup\theta_2\cup\theta_3 & 1 & 0 & 0.02 & 0.05  & 0.07\\
\hline
\end{array}
\end{equation*}

From the previous table of this numerical example, we see in column corresponding to $S_3(A)$ how the initial combined masses $m_{\mathcal{M}^f(\Theta)}((\theta_1\cup\theta_3)\cap\theta_2)\equiv S_1((\theta_1\cup\theta_3)\cap\theta_2)=0.05$, $m_{\mathcal{M}^f(\Theta)} (\theta_1\cap\theta_2\cap\theta_3) \equiv S_1(\theta_1\cap\theta_2\cap\theta_3)=0.16$, $m_{\mathcal{M}^f(\Theta)}(\theta_2\cap\theta_3)\equiv S_1(\theta_2\cap\theta_3)=0.19$ and 
$m_{\mathcal{M}^f(\Theta)}(\theta_1\cap\theta_2)\equiv S_1(\theta_1\cap\theta_2)=0.22$
are transferred (due to the constraint of $\mathcal{M}_3$) onto some elements of $D^\Theta$. 
We can easily check that the sum of the elements of the column for $S_3(A)$ is equal to $0.05+0.16+0.19+0.22=0.62$. \\

Because some elements of $D^\Theta$ are now equivalent due to the constraints of $\mathcal{M}_3$, we have to sum all the masses corresponding to same equivalent propositions. Thus after the final compression step, one gets the reduced hyper-power set $D^\Theta_{\mathcal{M}_3}$ having only 10 different elements with the following combined belief masses :
\begin{equation*}
\begin{array}{|l|r|}
\hline
\text{Element} \: A \: \text{of} \; D^\Theta_{\mathcal{M}_3} & m_{\mathcal{M}_3(\Theta)}(A)\\
\hline
\emptyset               & 0\\
\theta_1\cap\theta_3 & 0.12+0.03+0.02+0=0.17\\ 
\theta_3 & 0.16+0.07=0.23\\      
\theta_2 & 0.12\\  
(\theta_1\cap\theta_3)\cup\theta_2 & 0.01\\
\theta_2\cup\theta_3 & 0.05\\
\theta_1 & 0.12+0.04=0.16\\
\theta_1\cup\theta_3 & 0.08\\
\theta_1\cup\theta_2 & 0.11\\
\theta_1\cup\theta_2\cup\theta_3 &  0.07\\
\hline
\end{array}
\qquad \text{and} \qquad
\mathbf{D}_{\mathcal{M}_3}=\begin{bmatrix}
    0  &   0       &  0   &  0     \\    
     0  &   0    &   0  &   1     \\   
     0   &  0    &   1  &   1     \\
     0  &   1   &    0  &   0     \\
     0  &   1    &   0  &   1     \\
     0  &   1    &   1   &  1    \\
     1  &   0    &   0  &   1   \\
      1  &   0    &   1  &   1   \\
     1   &  1    &   0   &  1   \\
     1  &   1    &   1  &   1    
\end{bmatrix}
\end{equation*}
The basis vector $\mathbf{u}_{\mathcal{M}_3}$ is given by $\mathbf{u}_{\mathcal{M}_3}=[<1>    <2>    <3>    <13>  ]'$ and the encoding matrix $\mathbf{D}_{\mathcal{M}_3}$ is explicated just above.

\clearpage
\newpage
\subsubsection{Application of the DSm Hybrid rule on example 4 (Shafer's model)}
%---------------------------------------------------------------------------------

Here is the numerical result corresponding to example 4 with the hybrid-model $\mathcal{M}_4$ including all possible exclusivity constraints. This DSm hybrid model corresponds actually to the Shafer's model. One gets now
\vspace{0.5cm}

\begin{equation*}
\begin{array}{|l|rrrrr|r|}
\hline
\text{Element} \: A \: \text{of} \; D^\Theta&   \phi(A) & S_1(A) & S_2(A)& S_3(A) & m_{\mathcal{M}_4(\Theta)}(A)\\
\hline
\emptyset               & 0 & 0 & 0 & 0 & 0\\
\theta_1\cap\theta_2\cap\theta_3  \overset{\mathcal{M}_4}{\equiv}\emptyset & 0 & 0.16 & 0 & 0 & 0\\
\theta_2\cap\theta_3\overset{\mathcal{M}_4}{\equiv}\emptyset & 0 & 0.19 & 0 & 0 & 0\\
\theta_1\cap\theta_3\overset{\mathcal{M}_4}{\equiv}\emptyset & 0 & 0.12 & 0 & 0 & 0\\ 
(\theta_1\cup\theta_2)\cap\theta_3  \overset{\mathcal{M}_4}{\equiv} \emptyset & 0 & 0.01 & 0 & 0 .02 & 0\\ 
\theta_3 & 1 & 0.10 & 0 & 0.07  & 0.17\\      
\theta_1\cap\theta_2\overset{\mathcal{M}_4}{\equiv}\emptyset   & 0 & 0.22 & 0 & 0.02  & 0\\
(\theta_1\cup\theta_3)\cap\theta_2 \overset{\mathcal{M}_4}{\equiv}\emptyset & 0 & 0.05 & 0 & 0.02  & 0\\
(\theta_2\cup\theta_3)\cap\theta_1\overset{\mathcal{M}_4}{\equiv} \emptyset & 0 & 0 & 0 & 0.02  & 0\\
\{(\theta_1\cap\theta_2)\cup\theta_3\} \cap(\theta_1\cup\theta_2) \overset{\mathcal{M}_4}{\equiv} \emptyset
& 0 & 0  & 0 & 0  & 0\\
(\theta_1\cap\theta_2)\cup\theta_3 \overset{\mathcal{M}_4}{\equiv}  \theta_3 & 1 & 0  & 0 & 0.07  & 0.07\\  
\theta_2 & 1 & 0.03  & 0 & 0.09  & 0.12\\  
(\theta_1\cap\theta_3)\cup\theta_2 \overset{\mathcal{M}_4}{\equiv}  \theta_2 & 1 & 0  & 0 & 0.01  & 0.01\\
\theta_2\cup\theta_3 & 1 & 0  & 0 & 0.05  & 0.05\\
\theta_1 & 1 & 0.08  & 0 & 0.06  & 0.14\\
(\theta_2\cap\theta_3)\cup\theta_1 \overset{\mathcal{M}_4}{\equiv}  \theta_1 & 1 & 0.02  & 0 & 0.02  & 0.04\\
\theta_1\cup\theta_3 & 1 & 0.02  & 0 & 0.15 & 0.17\\
\theta_1\cup\theta_2 & 1 & 0 & 0.02 & 0.09  & 0.11\\
\theta_1\cup\theta_2\cup\theta_3 & 1 & 0 & 0.06 & 0.06  & 0.12\\
\hline
\end{array}
\end{equation*}
\vspace{0.5cm}
 
From the previous table of this numerical example, we see in column corresponding to $S_3(A)$ how the initial combined masses of the eight elements forced to the empty set by the constraints of the model $\mathcal{M}_4$ are transferred onto some elements of $D^\Theta$. 
We can easily check that the sum of the elements of the column for $S_3(A)$ is equal to $0.16+0.19+0.12+0.01+0.22+0.05+0=0.75$. \\

After the final compression step (i.e. the clustering of all equivalent propositions), one gets the reduced hyper-power set $D^\Theta_{\mathcal{M}_4}$ having only $2^3=8$ (corresponding to the classical power set $2^\Theta$) with the following combined belief masses:

\vspace{0.5cm}
\begin{equation*}
\begin{array}{|l|r|}
\hline
\text{Element} \: A \: \text{of} \; D^\Theta_{\mathcal{M}_4} & m_{\mathcal{M}_4(\Theta)}(A)\\
\hline
\emptyset               & 0\\
\theta_3 & 0.17+0.07=0.24\\      
\theta_2 & 0.12+0.01=0.13\\  
\theta_2\cup\theta_3 & 0.05\\
\theta_1 & 0.14+0.04=0.18\\
\theta_1\cup\theta_3 & 0.17\\
\theta_1\cup\theta_2 & 0.11\\
\theta_1\cup\theta_2\cup\theta_3 &  0.12\\
\hline
\end{array}
\qquad \text{and} \qquad
\mathbf{D}_{\mathcal{M}_4}=\begin{bmatrix}
    0  &   0       &  0    \\    
     0   &  0    &   1   \\
     0  &   1   &    0   \\
      0  &   1    &   1    \\
     1  &   0    &   0   \\
      1  &   0    &   1   \\
     1   &  1    &   0    \\
     1  &   1    &   1  
\end{bmatrix}
\end{equation*}

\vspace{0.5cm}

The basis vector $\mathbf{u}_{\mathcal{M}_4}$ is given by $\mathbf{u}_{\mathcal{M}_4}=[<1>    <2>    <3> ]'$ and the encoding matrix $\mathbf{D}_{\mathcal{M}_4}$ is explicated just above.

\clearpage
\newpage

\subsubsection{Application of the DSm Hybrid rule on example 5}
%---------------------------------------------------------------------------------

Here is the numerical result corresponding to example 5 with the hybrid-model $\mathcal{M}_5$ including the non-existential constraint $\theta_1\overset{\mathcal{M}_5}{\equiv}\emptyset$.
This non-existential constraint implies $\theta_1\cap\theta_2\cap\theta_3 \overset{\mathcal{M}_5}{\equiv}\emptyset$, $\theta_1\cap\theta_2 \overset{\mathcal{M}_5}{\equiv}\emptyset$, $\theta_1\cap\theta_3 \overset{\mathcal{M}_5}{\equiv}\emptyset$ and $(\theta_2\cup\theta_3)\cap\theta_1\overset{\mathcal{M}_5}{\equiv}\emptyset $.
One gets now with applying the DSm hybrid rule of combination:

 \vspace{0.5cm}

\begin{equation*}
\begin{array}{|l|rrrrr|r|}
\hline
\text{Element} \: A \: \text{of} \; D^\Theta&   \phi(A) & S_1(A) & S_2(A)& S_3(A) & m_{\mathcal{M}_5(\Theta)}(A)\\
\hline
\emptyset               & 0 & 0 & 0 & 0 & 0\\
%
% ii=1 ==> atom Aii=0  0  0  0  0  0  1
%Phi(A)=0 mf(A)=0.16 omega(A)=0 Sum2(A)=0.12 Sum3(A)=0 mai=0
\theta_1\cap\theta_2\cap\theta_3  \overset{\mathcal{M}_5}{\equiv}\emptyset & 0 & 0.16 & 0 & 0 & 0\\

% ii=2 ==> atom Aii=0  0  0  0  0  1  1
%Phi(A)=1 mf(A)=0.19 omega(A)=0 Sum2(A)=0.12 Sum3(A)=0 mai=0.19
\theta_2\cap\theta_3 & 1 & 0.19 & 0 & 0 & 0.19\\

% ii=3 ==> atom Aii=0  0  0  0  1  0  1
%Phi(A)=0 mf(A)=0.12 omega(A)=0 Sum2(A)=0.12 Sum3(A)=0 mai=0
\theta_1\cap\theta_3\overset{\mathcal{M}_5}{\equiv}\emptyset & 0 & 0.12  & 0 & 0 & 0\\ 

% ii=4 ==> atom Aii=0  0  0  0  1  1  1
%Phi(A)=1 mf(A)=0.01 omega(A)=0 Sum2(A)=0.12 Sum3(A)=0.02 mai=0.03
(\theta_1\cup\theta_2)\cap\theta_3  \overset{\mathcal{M}_5}{\equiv} \theta_2\cap\theta_3& 1 & 0.01  & 0 & 0 .02 & 0.03\\ 

% ii=5 ==> atom Aii=0  0  0  1  1  1  1
%Phi(A)=1 mf(A)=0.1 omega(A)=0 Sum2(A)=0.12 Sum3(A)=0.01 mai=0.11
\theta_3 & 1 & 0.10 & 0 & 0.01  & 0.11\\      

% ii=6 ==> atom Aii=0  0  1  0  0  0  1
%Phi(A)=0 mf(A)=0.22 omega(A)=0 Sum2(A)=0.12 Sum3(A)=0.02 mai=0
\theta_1\cap\theta_2\overset{\mathcal{M}_5}{\equiv}\emptyset   & 0 & 0.22 & 0 & 0.02  & 0\\

% ii=7 ==> atom Aii=0  0  1  0  0  1  1
%Phi(A)=1 mf(A)=0.05 omega(A)=0 Sum2(A)=0.12 Sum3(A)=0.02 mai=0.07
(\theta_1\cup\theta_3)\cap\theta_2 \overset{\mathcal{M}_5}{\equiv}\theta_2\cap\theta_3 & 1 & 0.05  & 0 & 0.02  & 0.07\\

% ii=8 ==> atom Aii=0  0  1  0  1  0  1
%Phi(A)=0 mf(A)=0 omega(A)=0 Sum2(A)=0.12 Sum3(A)=0.02 mai=0
(\theta_2\cup\theta_3)\cap\theta_1\overset{\mathcal{M}_5}{\equiv} \emptyset & 0 & 0 & 0 & 0.02  & 0\\

% ii=9 ==> atom Aii=0  0  1  0  1  1  1
%Phi(A)=1 mf(A)=0 omega(A)=0 Sum2(A)=0.12 Sum3(A)=0 mai=0
\{(\theta_1\cap\theta_2)\cup\theta_3\} \cap(\theta_1\cup\theta_2) \overset{\mathcal{M}_5}{\equiv} \theta_2\cap\theta_3
& 1 & 0 & 0 & 0  & 0\\

% ii=10 ==> atom Aii=0  0  1  1  1  1  1
%Phi(A)=1 mf(A)=0 omega(A)=0 Sum2(A)=0.12 Sum3(A)=0.07 mai=0.07
(\theta_1\cap\theta_2)\cup\theta_3 \overset{\mathcal{M}_5}{\equiv}  \theta_3 & 1 & 0 & 0 & 0.07  & 0.07\\  

% ii=11 ==> atom Aii=0  1  1  0  0  1  1
%Phi(A)=1 mf(A)=0.03 omega(A)=0 Sum2(A)=0.12 Sum3(A)=0.05 mai=0.08
\theta_2 & 1 & 0.03 & 0 & 0.05  & 0.08\\  
           
% ii=12 ==> atom Aii=0  1  1  0  1  1  1
%Phi(A)=1 mf(A)=0 omega(A)=0 Sum2(A)=0.12 Sum3(A)=0.01 mai=0.01
(\theta_1\cap\theta_3)\cup\theta_2 \overset{\mathcal{M}_5}{\equiv}  \theta_2 & 1 & 0 & 0 & 0.01  & 0.01\\

% ii=13 ==> atom Aii=0  1  1  1  1  1  1
%Phi(A)=1 mf(A)=0 omega(A)=0 Sum2(A)=0.12 Sum3(A)=0 mai=0
\theta_2\cup\theta_3 & 1 & 0  & 0 & 0  & 0\\

% ii=14 ==> atom Aii=1  0  1  0  1  0  1
%Phi(A)=0 mf(A)=0.08 omega(A)=0 Sum2(A)=0.12 Sum3(A)=0.08 mai=0
\theta_1 \overset{\mathcal{M}_5}{\equiv} \emptyset & 0 & 0.08 & 0.02 & 0.08  & 0\\

% ii=15 ==> atom Aii=1  0  1  0  1  1  1
%Phi(A)=1 mf(A)=0.02 omega(A)=0 Sum2(A)=0.12 Sum3(A)=0.02 mai=0.04
(\theta_2\cap\theta_3)\cup\theta_1 \overset{\mathcal{M}_5}{\equiv}  \theta_2\cap\theta_3 & 1 & 0.02 & 0 & 0.02  & 0.04\\

% ii=16 ==> atom Aii=1  0  1  1  1  1  1
%Phi(A)=1 mf(A)=0.02 omega(A)=0 Sum2(A)=0.12 Sum3(A)=0.17 mai=0.19
\theta_1\cup\theta_3 \overset{\mathcal{M}_5}{\equiv}  \theta_3 & 1 & 0.02 & 0.02 & 0.17 & 0.21\\

% ii=17 ==> atom Aii=1  1  1  0  1  1  1
%Phi(A)=1 mf(A)=0 omega(A)=0 Sum2(A)=0.12 Sum3(A)=0.09 mai=0.09
\theta_1\cup\theta_2 \overset{\mathcal{M}_5}{\equiv}  \theta_2 & 1 & 0 & 0.06 & 0.09  & 0.15\\

% ii=18 ==> atom Aii=1  1  1  1  1  1  1
%Phi(A)=1 mf(A)=0 omega(A)=1 Sum2(A)=0.12 Sum3(A)=0 mai=0.12
\theta_1\cup\theta_2\cup\theta_3 \overset{\mathcal{M}_5}{\equiv}  \theta_2 \cup \theta_3 & 1 & 0 & 0.04 & 0  & 0.04\\
\hline
\end{array}
\end{equation*}
\vspace{0.5cm}

From the previous table of this numerical example, we see in column corresponding to $S_3(A)$ how the initial combined masses of the 5 elements forced to the empty set by the constraints of the model $\mathcal{M}_5$ are transferred onto some elements of $D^\Theta$. 
We can easily check that the sum of the elements of the column for $S_3(A)$ is equal to $0+0.16+0.12+0.22+0+0.08=0.58$ (sum of  S1(A) for which $\phi(A)=0$). \\

After the final compression step (i.e. the clustering of all equivalent propositions), one gets the reduced hyper-power set $D^\Theta_{\mathcal{M}_5}$ having only 5 different elements according to:
\vspace{0.5cm}
\begin{equation*}
\begin{array}{|l|r|}
\hline
\text{Element} \: A \: \text{of} \; D^\Theta_{\mathcal{M}_5} & m_{\mathcal{M}_5(\Theta)}(A)\\
\hline
\emptyset               & 0\\
\theta_2\cap\theta_3 & 0.19+0.03+0.07+0+0.04=0.33\\      
\theta_3 & 0.11+0.07+0.21=0.39\\  
\theta_2 & 0.08+0.01+0.15=0.24\\
\theta_2\cup\theta_3 &  0+0.04=0.04\\
\hline
\end{array}
\qquad \text{and} \qquad
\mathbf{D}_{\mathcal{M}_5}=\begin{bmatrix}
      0  &   0   &  0    \\    
      0  &   0   &   1   \\
      0  &   1   &   1    \\
      1  &   0   &   1   \\
      1  &   1   &   1  
\end{bmatrix}
\end{equation*}
\vspace{0.5cm}

The basis vector $\mathbf{u}_{\mathcal{M}_5}$ is given by $\mathbf{u}_{\mathcal{M}_5}=[<2>    <3>    <23> ]'$. and the encoding matrix $\mathbf{D}_{\mathcal{M}_5}$ is explicated just above.

\clearpage
\newpage

\subsubsection{Application of the DSm Hybrid rule on example 6}
%---------------------------------------------------------------------------------
Here is the numerical result corresponding to example 6 with the hybrid-model $\mathcal{M}_6$ including the two non-existential constraint $\theta_1\overset{\mathcal{M}_6}{\equiv}\emptyset$ and $\theta_2\overset{\mathcal{M}_6}{\equiv}\emptyset$. This is a degenerate example actually, since no uncertainty arises in such trivial model. We just want to show here that the DSm hybrid rule still works in this example and provide a legitimist result.
By applying the DSm hybrid rule of combination, one now gets:
\vspace{0.5cm}
\begin{equation*}
\begin{array}{|l|rrrrr|r|}
\hline
\text{Element} \: A \: \text{of} \; D^\Theta&   \phi(A) & S_1(A) & S_2(A)& S_3(A) & m_{\mathcal{M}_6(\Theta)}(A)\\
\hline
\emptyset               & 0 & 0 & 0 & 0 & 0\\
%
% ii=1 ==> atom Aii=0  0  0  0  0  0  1
%Phi(A)=0 mf(A)=0.16 omega(A)=0 Sum2(A)=0.42 Sum3(A)=0 mai=0
\theta_1\cap\theta_2\cap\theta_3  \overset{\mathcal{M}_6}{\equiv}\emptyset & 0 & 0.16  & 0 & 0 & 0\\

% ii=2 ==> atom Aii=0  0  0  0  0  1  1
%Phi(A)=0 mf(A)=0.19 omega(A)=0 Sum2(A)=0.42 Sum3(A)=0 mai=0
\theta_2\cap\theta_3\overset{\mathcal{M}_6}{\equiv}\emptyset & 0 & 0.19 & 0 & 0 & 0\\

% ii=3 ==> atom Aii=0  0  0  0  1  0  1
%Phi(A)=0 mf(A)=0.12 omega(A)=0 Sum2(A)=0.42 Sum3(A)=0 mai=0
\theta_1\cap\theta_3\overset{\mathcal{M}_6}{\equiv}\emptyset & 0 & 0.12 & 0 & 0 & 0\\ 

% ii=4 ==> atom Aii=0  0  0  0  1  1  1
%Phi(A)=0 mf(A)=0.01 omega(A)=0 Sum2(A)=0.42 Sum3(A)=0.02 mai=0
(\theta_1\cup\theta_2)\cap\theta_3  \overset{\mathcal{M}_6}{\equiv} \emptyset & 0 & 0.01 & 0 & 0.02 & 0\\ 

% ii=5 ==> atom Aii=0  0  0  1  1  1  1
%Phi(A)=1 mf(A)=0.1 omega(A)=0 Sum2(A)=0.42 Sum3(A)=0.07 mai=0.17
\theta_3 & 1 & 0.10 & 0 & 0.07  & 0.17\\      

% ii=6 ==> atom Aii=0  0  1  0  0  0  1
%Phi(A)=0 mf(A)=0.22 omega(A)=0 Sum2(A)=0.42 Sum3(A)=0.02 mai=0
\theta_1\cap\theta_2\overset{\mathcal{M}_6}{\equiv}\emptyset   & 0 & 0.22 & 0 & 0.02  & 0\\

% ii=7 ==> atom Aii=0  0  1  0  0  1  1
%Phi(A)=0 mf(A)=0.05 omega(A)=0 Sum2(A)=0.42 Sum3(A)=0.02 mai=0
(\theta_1\cup\theta_3)\cap\theta_2 \overset{\mathcal{M}_6}{\equiv}\emptyset  & 0 & 0.05 & 0 & 0.02  & 0\\

% ii=8 ==> atom Aii=0  0  1  0  1  0  1
%Phi(A)=0 mf(A)=0 omega(A)=0 Sum2(A)=0.42 Sum3(A)=0.02 mai=0
(\theta_2\cup\theta_3)\cap\theta_1\overset{\mathcal{M}_6}{\equiv} \emptyset & 0 & 0 & 0 & 0.02  & 0\\

% ii=9 ==> atom Aii=0  0  1  0  1  1  1
%Phi(A)=0 mf(A)=0 omega(A)=0 Sum2(A)=0.42 Sum3(A)=0 mai=0
\{(\theta_1\cap\theta_2)\cup\theta_3\} \cap(\theta_1\cup\theta_2) \overset{\mathcal{M}_6}{\equiv} \emptyset
& 0 & 0 & 0 & 0  & 0\\

% ii=10 ==> atom Aii=0  0  1  1  1  1  1
%Phi(A)=1 mf(A)=0 omega(A)=0 Sum2(A)=0.42 Sum3(A)=0.07 mai=0.07
(\theta_1\cap\theta_2)\cup\theta_3 \overset{\mathcal{M}_6}{\equiv}  \theta_3 & 1 & 0 & 0 & 0.07  & 0.07\\  

% ii=11 ==> atom Aii=0  1  1  0  0  1  1
%Phi(A)=0 mf(A)=0.03 omega(A)=0 Sum2(A)=0.42 Sum3(A)=0.11 mai=0
\theta_2\overset{\mathcal{M}_6}{\equiv}\emptyset & 0 & 0.03 & 0.02 & 0.11  & 0\\  
           
% ii=12 ==> atom Aii=0  1  1  0  1  1  1
%Phi(A)=0 mf(A)=0 omega(A)=0 Sum2(A)=0.42 Sum3(A)=0.01 mai=0
(\theta_1\cap\theta_3)\cup\theta_2 \overset{\mathcal{M}_6}{\equiv} \emptyset & 0 & 0 & 0 & 0.01  & 0\\

% ii=13 ==> atom Aii=0  1  1  1  1  1  1
%Phi(A)=1 mf(A)=0 omega(A)=0 Sum2(A)=0.42 Sum3(A)=0.05 mai=0.05
\theta_2\cup\theta_3\overset{\mathcal{M}_6}{\equiv}  \theta_3 & 1 & 0 & 0.04 & 0.05 & 0.09\\

% ii=14 ==> atom Aii=1  0  1  0  1  0  1
%Phi(A)=0 mf(A)=0.08 omega(A)=0 Sum2(A)=0.42 Sum3(A)=0.08 mai=0
\theta_1 \overset{\mathcal{M}_6}{\equiv} \emptyset & 0 & 0.08 & 0 & 0.08  & 0\\

% ii=15 ==> atom Aii=1  0  1  0  1  1  1
%Phi(A)=0 mf(A)=0.02 omega(A)=0 Sum2(A)=0.42 Sum3(A)=0.02 mai=0
(\theta_2\cap\theta_3)\cup\theta_1 \overset{\mathcal{M}_6}{\equiv} \emptyset & 0 & 0.02 & 0 & 0.02  & 0\\

% ii=16 ==> atom Aii=1  0  1  1  1  1  1
%Phi(A)=1 mf(A)=0.02 omega(A)=0 Sum2(A)=0.42 Sum3(A)=0.19 mai=0.21
\theta_1\cup\theta_3 \overset{\mathcal{M}_6}{\equiv}  \theta_3 & 1 & 0.02 & 0.02 & 0.19 & 0.23\\

% ii=17 ==> atom Aii=1  1  1  0  1  1  1
%Phi(A)=0 mf(A)=0 omega(A)=0 Sum2(A)=0.42 Sum3(A)=0.12 mai=0
\theta_1\cup\theta_2 \overset{\mathcal{M}_6}{\equiv} \emptyset & 0 & 0 & 0.21 & 0.12  & 0\\

% ii=18 ==> atom Aii=1  1  1  1  1  1  1
%Phi(A)=1 mf(A)=0 omega(A)=1 Sum2(A)=0.42 Sum3(A)=0.08 mai=0.5
\theta_1\cup\theta_2\cup\theta_3 \overset{\mathcal{M}_6}{\equiv}  \theta_3 & 1 & 0 & 0.36 & 0.08  & 0.44\\
\hline
\end{array}
\end{equation*}
\vspace{0.5cm}

After the clustering of all equivalent propositions, one gets the reduced hyper-power set $D^\Theta_{\mathcal{M}_6}$ having only 2 different elements according to:
\vspace{0.5cm}
\begin{equation*}
\begin{array}{|l|r|}
\hline
\text{Element} \: A \: \text{of} \; D^\Theta_{\mathcal{M}_6} & m_{\mathcal{M}_6(\Theta)}(A)\\
\hline
\emptyset               & 0\\
\theta_3 & 0.17+0.07+0.09+0.23+0.44=1\\  
\hline
\end{array}
\end{equation*}
\vspace{0.5cm}
The encoding matrix $\mathbf{D}_{\mathcal{M}_6}$ and the basis vector $\mathbf{u}_{\mathcal{M}_6}$ for the elements of $D^\Theta_{\mathcal{M}_6}$ reduce to
$\mathbf{D}_{\mathcal{M}_6}=[0 1]'$ and $\mathbf{u}_{\mathcal{M}_6}=[<3>]$.

\clearpage
\newpage

\subsubsection{Application of the DSm Hybrid rule on example 7}
%---------------------------------------------------------------------------------

Here is the numerical result corresponding to example 7 with the hybrid-model $\mathcal{M}_7$ including the  {\it{mixed exclusivity and  non-existential}} constraint 
$( \theta_1\cap \theta_2)\cup\theta_3\overset{\mathcal{M}_7}{\equiv}\emptyset$. This mixed constraint implies
$\theta_1\cap\theta_2\cap\theta_3 \overset{\mathcal{M}_7}{\equiv}\emptyset$,
$\theta_1\cap\theta_2 \overset{\mathcal{M}_7}{\equiv}\emptyset$,
$\theta_1\cap\theta_3 \overset{\mathcal{M}_7}{\equiv}\emptyset$,
$\theta_2\cap\theta_3 \overset{\mathcal{M}_7}{\equiv}\emptyset$,
$(\theta_1\cup\theta_2)\cap\theta_3\overset{\mathcal{M}_7}{\equiv}\emptyset$,
$(\theta_1\cup\theta_3)\cap\theta_2\overset{\mathcal{M}_7}{\equiv}\emptyset$,
$(\theta_2\cup\theta_3)\cap\theta_1\overset{\mathcal{M}_7}{\equiv}\emptyset $,
$ \{(\theta_1\cap\theta_2)\cup\theta_3\} \cap(\theta_1\cup\theta_2)\overset{\mathcal{M}_7}{\equiv}\emptyset$ and
$\theta_3\overset{\mathcal{M}_7}{\equiv}\emptyset$.
By applying the DSm hybrid rule of combination, one now gets:
\vspace{0.5cm}
\begin{equation*}
\begin{array}{|l|rrrrr|r|}
\hline
\text{Element} \: A \: \text{of} \; D^\Theta&   \phi(A) & S_1(A) &  S_2(A)& S_3(A) & m_{\mathcal{M}_7(\Theta)}(A)\\
\hline
\emptyset               & 0 & 0 & 0 & 0 & 0\\
%
% ii=1 ==> atom Aii=0  0  0  0  0  0  1
%Phi(A)=0 mf(A)=0.16 omega(A)=0 Sum2(A)=0.25 Sum3(A)=0 mai=0
\theta_1\cap\theta_2\cap\theta_3  \overset{\mathcal{M}_7}{\equiv}\emptyset & 0 & 0.16  & 0 & 0 & 0\\

%% ii=2 ==> atom Aii=0  0  0  0  0  1  1
%%Phi(A)=0 mf(A)=0.19 omega(A)=0 Sum2(A)=0.25 Sum3(A)=0 mai=0
\theta_2\cap\theta_3\overset{\mathcal{M}_7}{\equiv}\emptyset & 0 & 0.19  & 0 & 0 & 0\\

%% ii=3 ==> atom Aii=0  0  0  0  1  0  1
%%Phi(A)=0 mf(A)=0.12 omega(A)=0 Sum2(A)=0.25 Sum3(A)=0 mai=0
\theta_1\cap\theta_3\overset{\mathcal{M}_7}{\equiv}\emptyset & 0 & 0.12  & 0 & 0 & 0\\ 

%% ii=4 ==> atom Aii=0  0  0  0  1  1  1
%%Phi(A)=0 mf(A)=0.01 omega(A)=0 Sum2(A)=0.25 Sum3(A)=0.02 mai=0
(\theta_1\cup\theta_2)\cap\theta_3  \overset{\mathcal{M}_7}{\equiv} \emptyset & 0 & 0.01  & 0 & 0 .02 & 0\\ 

%% ii=5 ==> atom Aii=0  0  0  1  1  1  1
%%Phi(A)=0 mf(A)=0.1 omega(A)=0 Sum2(A)=0.25 Sum3(A)=0.1 mai=0
\theta_3\overset{\mathcal{M}_7}{\equiv} \emptyset  & 0 & 0.10  & 0.03 & 0.10  & 0\\      

%% ii=6 ==> atom Aii=0  0  1  0  0  0  1
%%Phi(A)=0 mf(A)=0.22 omega(A)=0 Sum2(A)=0.25 Sum3(A)=0.02 mai=0
\theta_1\cap\theta_2\overset{\mathcal{M}_7}{\equiv}\emptyset   & 0 & 0.22  & 0 & 0.02  & 0\\

% ii=7 ==> atom Aii=0  0  1  0  0  1  1
%Phi(A)=0 mf(A)=0.05 omega(A)=0 Sum2(A)=0.25 Sum3(A)=0.02 mai=0
(\theta_1\cup\theta_3)\cap\theta_2 \overset{\mathcal{M}_7}{\equiv}\emptyset  & 0 & 0.05  & 0 & 0.02  & 0\\

% ii=8 ==> atom Aii=0  0  1  0  1  0  1
%Phi(A)=0 mf(A)=0 omega(A)=0 Sum2(A)=0.25 Sum3(A)=0.02 mai=0
(\theta_2\cup\theta_3)\cap\theta_1\overset{\mathcal{M}_7}{\equiv} \emptyset & 0 & 0  & 0 & 0.02  & 0\\

% ii=9 ==> atom Aii=0  0  1  0  1  1  1
%Phi(A)=0 mf(A)=0 omega(A)=0 Sum2(A)=0.25 Sum3(A)=0 mai=0
\{(\theta_1\cap\theta_2)\cup\theta_3\} \cap(\theta_1\cup\theta_2) \overset{\mathcal{M}_7}{\equiv} \emptyset
& 0 & 0  & 0 & 0  & 0\\

% ii=10 ==> atom Aii=0  0  1  1  1  1  1
%Phi(A)=0 mf(A)=0 omega(A)=0 Sum2(A)=0.25 Sum3(A)=0.07 mai=0
(\theta_1\cap\theta_2)\cup\theta_3 \overset{\mathcal{M}_7}{\equiv}  \emptyset & 0 & 0  & 0 & 0.07  & 0\\  

% ii=11 ==> atom Aii=0  1  1  0  0  1  1
%Phi(A)=1 mf(A)=0.03 omega(A)=0 Sum2(A)=0.25 Sum3(A)=0.09 mai=0.12
\theta_2 & 1 & 0.03  & 0 & 0.09  & 0.12\\  
           
% ii=12 ==> atom Aii=0  1  1  0  1  1  1
%Phi(A)=1 mf(A)=0 omega(A)=0 Sum2(A)=0.25 Sum3(A)=0.01 mai=0.01
(\theta_1\cap\theta_3)\cup\theta_2 \overset{\mathcal{M}_7}{\equiv} \theta_2 & 1 & 0  & 0 & 0.01  & 0.01\\

% ii=13 ==> atom Aii=0  1  1  1  1  1  1
%Phi(A)=1 mf(A)=0 omega(A)=0 Sum2(A)=0.25 Sum3(A)=0.05 mai=0.05
\theta_2\cup\theta_3\overset{\mathcal{M}_7}{\equiv}  \theta_2 & 1 & 0  & 0.06 & 0.05 & 0.11\\

% ii=14 ==> atom Aii=1  0  1  0  1  0  1
%Phi(A)=1 mf(A)=0.08 omega(A)=0 Sum2(A)=0.25 Sum3(A)=0.06 mai=0.14
\theta_1  & 1 & 0.08  & 0 & 0.06  & 0.14\\

% ii=15 ==> atom Aii=1  0  1  0  1  1  1
%Phi(A)=1 mf(A)=0.02 omega(A)=0 Sum2(A)=0.25 Sum3(A)=0.02 mai=0.04
(\theta_2\cap\theta_3)\cup\theta_1 \overset{\mathcal{M}_7}{\equiv} \theta_1 & 1 & 0.02  & 0 & 0.02  & 0.04\\

% ii=16 ==> atom Aii=1  0  1  1  1  1  1
%Phi(A)=1 mf(A)=0.02 omega(A)=0 Sum2(A)=0.25 Sum3(A)=0.22 mai=0.24
\theta_1\cup\theta_3 \overset{\mathcal{M}_7}{\equiv}  \theta_1 & 1 & 0.02  & 0.01 & 0.22 & 0.25\\

% ii=17 ==> atom Aii=1  1  1  0  1  1  1
%Phi(A)=1 mf(A)=0 omega(A)=0 Sum2(A)=0.25 Sum3(A)=0.09 mai=0.09
\theta_1\cup\theta_2 & 1 & 0  & 0.02 & 0.09  & 0.11\\

% ii=18 ==> atom Aii=1  1  1  1  1  1  1
%Phi(A)=1 mf(A)=0 omega(A)=1 Sum2(A)=0.25 Sum3(A)=0.06 mai=0.31
\theta_1\cup\theta_2\cup\theta_3 \overset{\mathcal{M}_7}{\equiv} \theta_1\cup \theta_2 & 1 & 0 & 0.16 & 0.06  & 0.22\\
\hline
\end{array}
\end{equation*}

\vspace{0.5cm}

After the clustering of all equivalent propositions, one gets the reduced hyper-power set $D^\Theta_{\mathcal{M}_6}$ having only 4 different elements according to:
\vspace{0.5cm}

\begin{equation*}
\begin{array}{|l|r|}
\hline
\text{Element} \: A \: \text{of} \; D^\Theta_{\mathcal{M}_7} & m_{\mathcal{M}_7(\Theta)}(A)\\
\hline
\emptyset               & 0\\
\theta_2 & 0.12+0.01+0.11=0.24\\  
\theta_1 & 0.14+0.04+0.25=0.43\\  
\theta_1\cup \theta_2 & 0.11 + 0.22=0.33\\  
\hline
\end{array}
\end{equation*}

\vspace{0.5cm}

The basis vector $\mathbf{u}_{\mathcal{M}_7}$ and the encoding matrix $\mathbf{D}_{\mathcal{M}_7}$ for the elements of $D^\Theta_{\mathcal{M}_7}$ are given by  
\begin{equation*}
\mathbf{u}_{\mathcal{M}_7}=[<1>   <2> ]' \qquad \text{and}\qquad 
\mathbf{D}_{\mathcal{M}_7}=\begin{bmatrix}
      0  &   0    \\
      0  &   1      \\
      1  &   0     \\
      1  &   1  
\end{bmatrix}
\end{equation*}
 
\clearpage
\newpage

%----------------------------------------------------------------------------------------------------------------------------
\subsection{Example  with more general basic belief assignments $m_1(.)$ and $m_2(.)$}
%----------------------------------------------------------------------------------------------------------------------------

We present in this section the numerical results of the DSm hybrid rule of combination applied upon the seven previous models $\mathcal{M}_i$, $i=1,...,7$ with two general basic belief assignments $m_1(.)$ and $m_2(.)$ such that $m_1(A)>0$ and $m_2(A)>0$ for all $A\neq\emptyset\in D^{\Theta=\{\theta_1,\theta_2,\theta_3\}}$. We just provide here results. The verification is left to the reader. The following table presents the numerical values chosen for $m_1(.)$ and $m_2(.)$ and the result of the fusion obtained by the classical DSm rule of combination

\begin{equation*}
\begin{array}{|l|ll|l|}
\hline
\text{Element} \: A \: \text{of} \; D^\Theta &   m_1(A) & m_2(A) & m_{\mathcal{M}^f}(A) \\
\hline
\emptyset                                               &       0 &     0  &   0   \\
\theta_1\cap\theta_2\cap\theta_3     & 0.01 & 0.40 & 0.4389 \\
\theta_2\cap\theta_3                           & 0.04 &0.03  & 0.0410 \\
\theta_1\cap\theta_3                           & 0.03 & 0.04 & 0.0497 \\
(\theta_1\cup\theta_2)\cap\theta_3  & 0.01 & 0.02 & 0.0257 \\
\theta_3                                                 & 0.03 & 0.04 & 0.0311 \\                          
\theta_1\cap\theta_2                           & 0.02 & 0.20 & 0.1846 \\
(\theta_1\cup\theta_3)\cap\theta_2  & 0.02 & 0.01 & 0.0156 \\
(\theta_2\cup\theta_3)\cap\theta_1  & 0.03 & 0.04 & 0.0459 \\
\{(\theta_1\cap\theta_2)\cup\theta_3\} \cap(\theta_1\cup\theta_2) & 0.04 & 0.03 & 0.0384 \\
(\theta_1\cap\theta_2)\cup\theta_3  & 0.04 & 0.03 & 0.0296 \\
\theta_2                                                 & 0.02 & 0.01 & 0.0084\\                    
(\theta_1\cap\theta_3)\cup\theta_2  & 0.01 & 0.02 & 0.0221 \\
\theta_2\cup\theta_3                           & 0.20 & 0.02 & 0.0140 \\
\theta_1                                                 & 0.01 & 0.02 & 0.0109 \\
(\theta_2\cap\theta_3)\cup\theta_1 & 0.02 & 0.01 & 0.0090 \\
\theta_1\cup\theta_3                          & 0.04 & 0.03 & 0.0136 \\
\theta_1\cup\theta_2                          & 0.03 & 0.04 & 0.0175 \\
\theta_1\cup\theta_2\cup\theta_3   & 0.40 & 0.01 & 0.0040  \\
\hline
\end{array}
\end{equation*}

The following table shows the results obtained by the DSm hybrid rule of combination before the final compression step of all redundant propositions for the DSm hybrid models presented in the previous examples.

\begin{equation*}
\begin{array}{|l|lllllll|}
\hline
\text{Element} \: A \: \text{of} \; D^\Theta & m_{\mathcal{M}_1}(A) & m_{\mathcal{M}_2}(A) & m_{\mathcal{M}_3}(A) & m_{\mathcal{M}_4}(A) & m_{\mathcal{M}_5}(A) & m_{\mathcal{M}_6}(A) & m_{\mathcal{M}_7}(A)\\
\hline
\emptyset                                               & 0 & 0 & 0 & 0 & 0 & 0 & 0\\
\theta_1\cap\theta_2\cap\theta_3     & 0 & 0 & 0 & 0 & 0 & 0 & 0\\
\theta_2\cap\theta_3                           & 0.0573 & 0.0573 & 0 & 0 & 0.0573 & 0 & 0\\
\theta_1\cap\theta_3                           & 0.0621 & 0.0621 & 0.0621& 0 & 0 & 0 & 0\\
(\theta_1\cup\theta_2)\cap\theta_3  & 0.0324 & 0.0324 & 0.0335 & 0 & 0.0334 & 0 & 0\\
\theta_3                                                 & 0.0435 & 0.0435 & 0.0460 & 0.0494 & 0.0459 & 0.0494 & 0\\                          
\theta_1\cap\theta_2                           & 0.1946 & 0 & 0 & 0 & 0 & 0 & 0\\
(\theta_1\cup\theta_3)\cap\theta_2  & 0.0323 & 0.0365 & 0 & 0 & 0.0365 & 0 & 0\\
(\theta_2\cup\theta_3)\cap\theta_1  & 0.0651 & 0.0719 & 0.0719 & 0 & 0 & 0 & 0\\
\{(\theta_1\cap\theta_2)\cup\theta_3\} \cap(\theta_1\cup\theta_2) & 0.0607 & 0.0704 & 0.0743 & 0 & 0.0764 & 0 & 0\\
(\theta_1\cap\theta_2)\cup\theta_3  & 0.0527 & 0.0613 & 0.0658 & 0.0792 & 0.0687 & 0.0792 & 0\\
\theta_2                                                 & 0.0165 & 0.0207 & 0.0221 & 0.0221& 0.0207 & 0 & 0.0221\\                    
(\theta_1\cap\theta_3)\cup\theta_2  & 0.0274 & 0.0309 & 0.0340 & 0.0375 & 0.0329 & 0 & 0.0375\\
\theta_2\cup\theta_3                           & 0.0942 & 0.1346 & 0.1471 & 0.1774 & 0.1518 & 0.1850 & 0.1953\\
\theta_1                                                 & 0.0151 & 0.0175 & 0.0175 & 0.0195 & 0 & 0 & 0.0195\\
(\theta_2\cap\theta_3)\cup\theta_1 & 0.0182 & 0.0229 & 0.0243 & 0.0295 & 0.0271 & 0 & 0.0295\\
\theta_1\cup\theta_3                          & 0.0299 & 0.0385 & 0.0419 & 0.0558 & 0.0489 & 0.0589 & 0.0631\\
\theta_1\cup\theta_2                          & 0.0299 & 0.0412 & 0.0452 & 0.0544 & 0.0498 & 0 & 0.0544\\
\theta_1\cup\theta_2\cup\theta_3   & 0.1681 & 0.2583 & 0.3143 & 0.4752 & 0.3506 & 0.6275 & 0.5786\\
\hline
\end{array}
\end{equation*}

The next tables present the final result of the DSm hybrid rule of combination after the compression step (the merging of all equivalent redundant propositions) presented in previous examples.

\clearpage
\newpage
\twocolumn

\begin{equation*}
\begin{array}{|l|r|}
\hline
\text{Element} \: A \: \text{of} \; D^\Theta_{\mathcal{M}_1} & m_{\mathcal{M}_1(\Theta)}(A)\\
\hline
\emptyset                                               & 0 \\
\theta_2\cap\theta_3                           & 0.0573 \\
\theta_1\cap\theta_3                           & 0.0621 \\
(\theta_1\cup\theta_2)\cap\theta_3  & 0.0324 \\
\theta_3                                                 & 0.0435 \\                          
\theta_1\cap\theta_2                           & 0.1946 \\
(\theta_1\cup\theta_3)\cap\theta_2  & 0.0323 \\
(\theta_2\cup\theta_3)\cap\theta_1  & 0.0651 \\
\{(\theta_1\cap\theta_2)\cup\theta_3\} \cap(\theta_1\cup\theta_2) & 0.0607 \\
(\theta_1\cap\theta_2)\cup\theta_3  & 0.0527 \\
\theta_2                                                 & 0.0165 \\                    
(\theta_1\cap\theta_3)\cup\theta_2  & 0.0274 \\
\theta_2\cup\theta_3                           & 0.0942 \\
\theta_1                                                 & 0.0151 \\
(\theta_2\cap\theta_3)\cup\theta_1 & 0.0182 \\
\theta_1\cup\theta_3                          & 0.0299 \\
\theta_1\cup\theta_2                          & 0.0299 \\
\theta_1\cup\theta_2\cup\theta_3   & 0.1681 \\
\hline
\end{array}
\end{equation*}

\vspace{0.5cm}

\begin{equation*}
\begin{array}{|l|r|}
\hline
\text{Element} \: A \: \text{of} \; D^\Theta_{\mathcal{M}_2} & m_{\mathcal{M}_2(\Theta)}(A)\\
\hline
\emptyset               & 0\\
\theta_2\cap\theta_3 &  0.0938\\
\theta_1\cap\theta_3 & 0.1340\\ 
(\theta_1\cup\theta_2)\cap\theta_3  & 0.1028\\ 
\theta_3 & 0.1048\\      
\theta_2 & 0.0207\\  
(\theta_1\cap\theta_3)\cup\theta_2 & 0.0309\\
\theta_2\cup\theta_3 & 0.1346\\
\theta_1 & 0.0175\\
(\theta_2\cap\theta_3)\cup\theta_1& 0.0229\\
\theta_1\cup\theta_3 & 0.0385\\
\theta_1\cup\theta_2 & 0.0412\\
\theta_1\cup\theta_2\cup\theta_3 &  0.2583\\
\hline
\end{array}
\end{equation*}

\vspace{0.5cm}

\begin{equation*}
\begin{array}{|l|r|}
\hline
\text{Element} \: A \: \text{of} \; D^\Theta_{\mathcal{M}_3} & m_{\mathcal{M}_3(\Theta)}(A)\\
\hline
\emptyset               & 0\\
\theta_1\cap\theta_3 & 0.2418\\ 
\theta_3 & 0.1118\\      
\theta_2 & 0.0221\\  
(\theta_1\cap\theta_3)\cup\theta_2 & 0.0340\\
\theta_2\cup\theta_3 & 0.1471\\
\theta_1 & 0.0418\\
\theta_1\cup\theta_3 & 0.0419\\
\theta_1\cup\theta_2 & 0.0452\\
\theta_1\cup\theta_2\cup\theta_3 &  0.3143\\
\hline
\end{array}
\end{equation*}

\vspace{0.5cm}

\begin{equation*}
\begin{array}{|l|r|}
\hline
\text{Element} \: A \: \text{of} \; D^\Theta_{\mathcal{M}_4} & m_{\mathcal{M}_4(\Theta)}(A)\\
\hline
\emptyset               & 0\\
\theta_3 & 0.1286\\      
\theta_2 & 0.0596\\  
\theta_2\cup\theta_3 & 0.1774\\
\theta_1 & 0.0490\\
\theta_1\cup\theta_3 & 0.0558\\
\theta_1\cup\theta_2 & 0.0544\\
\theta_1\cup\theta_2\cup\theta_3 & 0.4752\\
\hline
\end{array}
\end{equation*}

\vspace{0.5cm}

\begin{equation*}
\begin{array}{|l|r|}
\hline
\text{Element} \: A \: \text{of} \; D^\Theta_{\mathcal{M}_5} & m_{\mathcal{M}_5(\Theta)}(A)\\
\hline
\emptyset               & 0\\
\theta_2\cap\theta_3 & 0.2307\\      
\theta_3 & 0.1635\\  
\theta_2 & 0.1034\\
\theta_2\cup\theta_3 &  0.5024\\
\hline
\end{array}
\end{equation*}

\vspace{0.5cm}

\begin{equation*}
\begin{array}{|l|r|}
\hline
\text{Element} \: A \: \text{of} \; D^\Theta_{\mathcal{M}_6} & m_{\mathcal{M}_6(\Theta)}(A)\\
\hline
\emptyset               & 0\\
\theta_3 & 1\\  
\hline
\end{array}
\end{equation*}

\vspace{0.5cm}

\begin{equation*}
\begin{array}{|l|r|}
\hline
\text{Element} \: A \: \text{of} \; D^\Theta_{\mathcal{M}_7} & m_{\mathcal{M}_7(\Theta)}(A)\\
\hline
\emptyset               & 0\\
\theta_2 & 0.2549\\  
\theta_1 & 0.1121\\  
\theta_1\cup \theta_2 & 0.6330\\  
\hline
\end{array}
\end{equation*}
\clearpage
\newpage
\onecolumn

%***********************************************************************
\subsection{DSm hybrid rule versus Dempster's rule of combination}
%***********************************************************************

We discuss and compare here the DSm hybrid rule of combination with respect to the Dempster's rule of combination and its alternative proposed in the literature based on  the Dempster-Shafer Theory (DST) framework which is frequently adopted in many fusion/expert systems. It is necessary to first recall briefly the basis of the DST \cite{Shafer_1976}.\\

\subsubsection{Brief introduction to the DST}

The DST  starts by assuming an exhaustive and exclusive frame of discernment of the problem under consideration $\Theta=\{\theta_1,\theta_2,\ldots,\theta_n\}$. This corresponds to the Shafer's model of the problem. The Shafer's model is nothing more but the DSm model including all possible exclusivity constraints.  The Shafer's model assumes actually that an ultimate refinement of the problem is possible so that $\theta_i$ are well precisely defined/identified in such a way  that we are sure that they are exclusive and exhaustive. From this Shafer's model,  a basic belief assignment (bba) $m (.): 2^\Theta \rightarrow  [0, 1]$  associated to a given body of evidence $\mathcal{B}$ (also called sometimes corpus of evidence) is defined by
\begin{equation}
m(\emptyset)=0  \qquad \text{and}\qquad     \sum_{A\in 2^\Theta} m(A) = 1                            
\end{equation}
\noindent
where $2^\Theta$ is called the {\it{power set}} of $\Theta$, i.e. the set of all subsets of $\Theta$. The set of all propositions $A\in 2^\Theta$ such that $m(A)>0$ is called the core of $m(.)$ and is denoted $\mathcal{K}(m)$. From any bba, one defines the belief and plausibility functions of $A\subseteq\Theta$ as
\begin{equation}
\text{Bel}(A) = \sum_{B\in 2^\Theta, B\subseteq A} m(B)
\label{Belg}
\end{equation}
\begin{equation}
\text{Pl}(A) = \sum_{B\in 2^\Theta, B\cap A\neq\emptyset} m(B)=1- \text{Bel}(\bar{A})
\label{Plg}
\end{equation}

\subsubsection{The Dempster's rule of combination}

Now let $\text{Bel}_1(.)$ and $\text{Bel}_2(.)$ be two belief functions over the same frame of discernment $\Theta$ and their corresponding bba $m_1(.)$ and $m_2(.)$ provided by two distinct bodies of evidence $\mathcal{B}_1$ and $\mathcal{B}_2$. Then the combined global belief function denoted $\text{Bel}(.)= \text{Bel}_1(.)\oplus \text{Bel}_2(.)$ is obtained by combining the basic belief assignments (called also sometimes  {\it{information granules}} in the literature) $m_1(.)$ and $m_2(.)$ through the following Dempster's rule of combination $[m_{1}\oplus m_{2}](\emptyset)=0$ and $\forall B\neq\emptyset \in 2^\Theta$,
 \begin{equation}
[m_{1}\oplus m_{2}](B) = 
\frac{\sum_{X\cap Y=B}m_{1}(X)m_{2}(Y)}{1-\sum_{X\cap Y=\emptyset} m_{1}(X) m_{2}(Y)} 
\label{eq:DSR}
 \end{equation}
 
The notation $\sum_{X\cap Y=B}$ represents the sum over all $X, Y \in 2^\Theta$ such that $X\cap Y=B$. 
The orthogonal sum $m (.)\triangleq [m_{1}\oplus m_{2}](.)$ is considered as a basic belief assignment if and only if the denominator in equation \eqref{eq:DSR} is non-zero. The term $k_{12}\triangleq \sum_{X\cap Y=\emptyset} m_{1}(X) m_{2}(Y)$ is called degree of conflict between the sources $\mathcal{B}_1$ and $\mathcal{B}_2$. When $k_{12}=1$,  the orthogonal sum $m (.)$ does not exist and the bodies of evidences $\mathcal{B}_1$ and $\mathcal{B}_2$ are said to be in {\it{full contradiction}}. Such a case can arise when there exists $A \subset \Theta$ such that $\text{Bel}_1(A) =1$ and $\text{Bel}_2(\bar{A}) = 1$. Same kind of trouble can occur also with the {\it{Optimal Bayesian Fusion Rule}} (OBFR) \cite{Dezert_2001a,Dezert_2001b}.\\

The DST is attractive for the {\it{Data Fusion community}} because it gives a nice mathematical model for ignorance and it includes the Bayesian theory as a special case \cite{Shafer_1976} (p. 4). Although very appealing, the DST presents some weaknesses and limitations because of its model itself, the theoretical justification of the Dempster's rule of combination but also because of our confidence to trust the result of Dempster's rule of combination when the conflict becomes important between sources ($k_{12} \nearrow 1$).\\

\subsubsection{Alternatives of the Dempster's rule of combination in the DST framework}

The Dempster's rule of combination has however been {\it{a posteriori}} justified  by the Smet's axiomatic of the Transferable Belief Model (TBM) in \cite{Smets_1994}. But we must also emphasize here that an infinite number of possible rules of combinations can be built from the Shafer's model  following ideas initially proposed by Lef\`evre, Colot and Vannoorenberghe in \cite{Lefevre_2002} and corrected here as follows:
\begin{itemize}
\item one first has to compute $m(\emptyset)$ by
$$ m(\emptyset) \triangleq \sum_{A\cap B =\emptyset}m_1(A)m_2(B) $$
\item then one redistributes $m(\emptyset)$ on all $(A\neq\emptyset)\subseteq \Theta$ with some given coefficients $w_m(A)\in[0,1]$ such that $\sum_{A\subseteq \Theta} w_m(A)=1$ according to
\begin{equation}
\begin{cases}
w_m(\emptyset)m(\emptyset) \rightarrow m(\emptyset)\\
m(A) + w_m(A)m(\emptyset) \rightarrow m(A), \forall A\neq\emptyset
\end{cases}
\label{eq:CEV}
\end{equation}
\end{itemize}
The particular choice of the set of coefficients $w_m(.)$ provides a particular rule of combination.
Actually there exists an infinite number of possible rules of combination. Some rules can be better justified than others depending on their ability to or not to preserve the associativity and commutativity properties of the combination. It can be easily shown in  \cite{Lefevre_2002}  that such general procedure provides all existing rules developed in the literature from the Shafer's model as alternative to the primeval Dempster's rule of combination depending on the choice of coefficients $w(A)$. As examples:
\begin{itemize}
\item
 the Dempster's rule of combination can be obtained from \eqref{eq:CEV} by choosing \cite{Lefevre_2002} $\forall A\neq\emptyset$
$$w_m(\emptyset)=0\qquad\text{and}\qquad w_m(A)=m(A)/(1-m(\emptyset))$$ 
\item
the Yager's rule of combination is obtained by choosing \cite{Yager_1987,Lefevre_2002} $$w_m(\Theta)=1$$ 
\item
 the Smets' rule of combination \cite{Smets_1990,Lefevre_2002}  is obtained by accepting the possibility to deal with bba such that $m(\emptyset)>0$ and thus by choosing $$w_m(\emptyset)=1$$
\item
with the Lef\'evre and al. formalism \cite{Lefevre_2002} and when $m(\emptyset)>0$,  the Dubois and Prade's rule of combination \cite{Dubois_1998,Lefevre_2002}  is obtained by choosing
$$\forall A\subseteq \mathcal{P}, \qquad w_m(A)=\frac{\sum_{\overset{A_1,A_2\mid A_1\cup A_2 = A}{A_1\cap A_2=\emptyset}}m^{\star}}{m(\emptyset}$$
\noindent where $m^{\star}\triangleq m_1(A_1)m_2(A_2)$  corresponds to the partial conflicting mass which is assigned to $A_1\cup A_2$ and where $\mathcal{P}$ is the set of all subsets of $2^\Theta$ on which the conflicting mass is distributed defined by

$$\mathcal{P}\triangleq\{A\in 2^\Theta \mid \exists A_1\in \mathcal{K}(m_1),  \exists A_2\in \mathcal{K}(m_2), A_1\cup A_2=A \; \text{and} \; A_1\cap A_2=\emptyset \}$$

The computation of the weighting factors $w_m(A)$ of the Dubois and Prade's rule of combination does not depend only on propositions they are associated with, but also on belief mass functions which have cause the partial conflicts. Thus the belief mass functions leading to the conflict allow to compute that part of conflicting mass which must be assigned to the subsets of $\mathcal{P}$ \cite{Lefevre_2002}. The Yager's rule coincides with the Dubois and Prade's rule of combination when choosing $\mathcal{P}=\{\Theta\}$. 

\end{itemize}

\subsubsection{DSm hybrid rule is not equivalent to the Dempster's rule of combination}

In its essence, the DSm hybrid rule of combination is  close to the Dubois and Prade's rule of combination but more general and precise because it works on $D^\Theta  \supset 2^\Theta$ and allows us to include all possible exclusivity and non-existential constraints for the model one has to work with. The advantage of using the DSm hybrid rule is that it does not require the calculation of weighting factors neither the normalization. The DSm hybrid rule of combination is definitely not equivalent to the Dempster's rule of combination as one can easily prove in the following very simple example:\\

Let consider $\Theta=\{\theta_1,\theta_2\}$ and the two sources in full contradiction providing the following basic belief assignments
$$m_1(\theta_1)=1 \qquad m_1(\theta_2)=0$$
$$m_2(\theta_1)=0 \qquad m_2(\theta_2)=1$$

Using the classic DSm rule of combination working with the free DSm model $\mathcal{M}^f$, one gets 

$$m_{\mathcal{M}^f}(\theta_1)=0 \qquad m_{\mathcal{M}^f}(\theta_2)=0 \qquad m_{\mathcal{M}^f}(\theta_1 \cap \theta_2)=1 \qquad m_{\mathcal{M}^f}(\theta_1 \cup \theta_2)=0$$

If one forces $\theta_1$ and $\theta_2$ to be exclusive to work with the Shafer's model $\mathcal{M}^0$, then the Dempster's rule of combination can not be applied in this limit case because of the full contradiction of the two sources of information. One gets the undefined operation 0/0. But the DSm hybrid rule can be applied in such limit case because it transfers the mass of this empty set ($\theta_1 \cap \theta_2\equiv \emptyset$ because of the choice of the model $\mathcal{M}^0$) to non-empty set(s), and one gets:

$$m_{\mathcal{M}^0}(\theta_1)=0 \qquad m_{\mathcal{M}^0}(\theta_2)=0 \qquad m_{\mathcal{M}^0}(\theta_1 \cap \theta_2)=0 \qquad m_{\mathcal{M}^0}(\theta_1 \cup \theta_2)=1$$

\noindent
This result is coherent in this very simple case with the Yager's and Dubois-Prade's rule of combination.\\

Now let examine the behavior of the numerical result when introducing a small variation $\epsilon > 0$ on initial basic belief assignments $m_1(.)$ and $m_2(.)$ as follows:

$$m_1(\theta_1)=1-\epsilon \qquad m_1(\theta_2)=\epsilon\qquad \text{and}\qquad m_2(\theta_1)=\epsilon \qquad m_2(\theta_2)=1-\epsilon $$

As shown on figure 1, $\lim_{\epsilon\rightarrow 0} m_{DS}(.)$, where $m_{DS}(.)$ is the result obtained from the Dempster's rule of combination, is given by

$$m_{DS}(\theta_1)=0.5 \qquad m_{DS}(\theta_2)=0.5\qquad m_{DS}(\theta_1\cap\theta_2)=0 \qquad m_{DS}(\theta_1\cup\theta_2)=0$$

This result is very questionable because it assigns same belief on $\theta_1$ and $\theta_2$ which is more informational than to assign all the belief to the total ignorance. 
The assignment of the belief to the total ignorance appears to be more justified from our point of view because it  properly reflects the almost total contradiction between the two sources and in such cases, it seems legitimist that the information can be drawn from the fusion.
When we apply the  DSm hybrid rule of combination (using the Shafer's model $\mathcal{M}^0$), one gets the expected belief assignment on the total ignorance, i.e. $m_{\mathcal{M}^0}(\theta_1 \cup \theta_2)=1$. The figure below shows the evolution of bba on $\theta_1$, $\theta_2$ and $\theta_1\cup\theta_2$ with $\epsilon$ obtained with classical Dempster's rule and DSm hybrid rule based on Shafer's model $\mathcal{M}^0$ (i.e. $\theta_1 \cap \theta_2\overset{\mathcal{M}_0}{\equiv} \emptyset$) .

\begin{figure}[h]
\centering
\includegraphics[width=17cm,height=6cm]{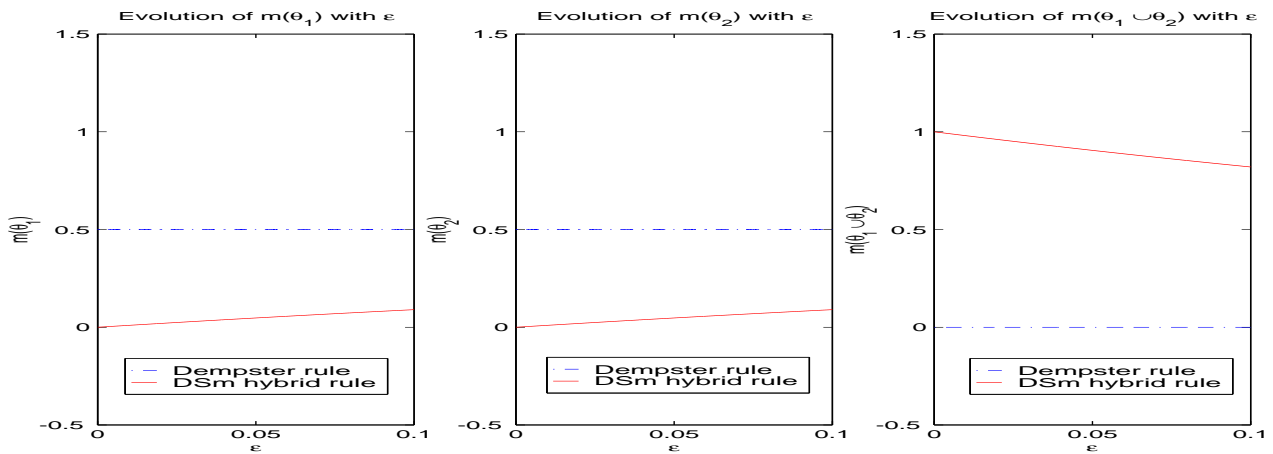}
 \caption{Comparison of Dempster's rule with the DSm hybrid rule on $\Theta=\{\theta_1,\theta_2\}$}
 \label{fig:2}
 \end{figure}

\clearpage
\newpage

%**************************
\section{Dynamic fusion}
%**************************

The DSm hybrid rule of combination presented in this paper has been developed for {\it{static}} problems/models, but is also directly applicable for easily handling {\it{dynamic fusion}} problems in real time as well, since at each temporal change of the models, one can still apply such hybrid rule. If $D^\Theta$ changes, due to the dynamicity of the frame $\Theta$, from time $t_l$ to time $t_{l+1}$, i.e. some of its elements which at time $t_l$ were not empty become (or are proved) empty at time $t_{l+1}$, or vice versa: if new elements, empty at time $t_l$, arise non-empty at time $t_{l+1}$, this DSm hybrid rule can be applied again at each change. If $\Theta$ tests the same but its set of focal (i.e. non-empty) elements of $D^\Theta$ increases, then again apply the
DSm hybrid rule.

\subsection{Example 1}

Let's consider the testimony fusion problem\footnote{This problem has been proposed to the authors in a private communication by L. Cholvy in 2002.} with the frame 
$$\Theta(t_l)\triangleq\{\theta_1\equiv\text{young},\theta_2\equiv\text{old},\theta_3\equiv\text{white hairs}\}$$
\noindent
with the following two basic belief assignments

$$m_1(\theta_1)=0.5 \qquad m_1(\theta_3)=0.5$$
$$m_2(\theta_2)=0.5 \qquad m_2(\theta_3)=0.5$$

\noindent
By applying the classical DSm fusion rule, one then gets

$$m_{\mathcal{M}^f(\Theta(t_l))}(\theta_1\cap\theta_2)=0.25 \qquad m_{\mathcal{M}^f(\Theta(t_l))}(\theta_1\cap\theta_3)=0.25 \qquad
m_{\mathcal{M}^f(\Theta(t_l))}(\theta_2\cap\theta_3)=0.25 \qquad
m_{\mathcal{M}^f(\Theta(t_l))}(\theta_3)=0.25$$

\noindent
Suppose now that at time $t_{l+1}$, one knows that young people don't have white hairs (i.e $\theta_1\cap\theta_3\equiv \emptyset$). How can we update the previous fusion result with this new information on the model of the problem ?
We solve it with the DSm hybrid rule, which transfers the mass of the empty sets (imposed by the constraints on the new model $\mathcal{M}$ available at time $t_{l+1}$) to the non-empty sets of $D^\Theta$, going on the track of the DSm classic rule. Using the DSm hybrid rule with the constraint $\theta_1\cap\theta_3\equiv \emptyset$, one then gets:

$$m_{\mathcal{M}}(\theta_1\cap\theta_2)=0.25 \qquad 
m_{\mathcal{M}}(\theta_2\cap\theta_3)=0.25 \qquad
m_{\mathcal{M}}(\theta_3)=0.25\qquad
$$
and the mass $m_{\mathcal{M}}(\theta_1\cap\theta_3)=0$, because $\theta_1\cap\theta_3=\{\text{young}\}\cap\{\text{white hairs}\}\overset{\mathcal{M}}{\equiv}\emptyset $ and its previous mass $m_{\mathcal{M}^f(\Theta(t_l))}(\theta_1\cap\theta_3)=0.25$ is transferred to $m_{\mathcal{M}}(\theta_1\cup\theta_3)=0.25$ by  the DSm hybrid rule.

%************************
\subsection{Example 2}
%************************

Let $\Theta(t_l)=\{\theta_1,\theta_2,\ldots,\theta_n\}$ be a list of suspects and let consider two observers who eyewitness the scene of plunder at a museum in Bagdad and who testify to the radio and TV the identities of thieves using the basic beliefs assignments $m_1(.)$ and $m_2(.)$ defined on $D^{\Theta(t_l)}$, where $t_l$ represents the time of the observation. Afterwards, at time $t_{l+1}$, one finds out that one suspect, among this list $\Theta(t_l)$, say $\theta_i$, could not be a suspect because he was on duty in another place, evidence which was certainly confirmed. Therefore he has to be taken off the suspect list $\Theta(t_l)$, and a new frame of discernment is resulting $\Theta(t_{l+1})$. If this one changes again, one applies again the DSm hybrid of combining of evidences, and so on. This is a typically dynamical example where models change with time and where one needs to adapt fusion results with current model over time. In the meantime, one can also take into account new observations/testimonies in the DSm hybrid fusion rule as soon as they become available to the fusion system. If $\Theta$ and $D^\Theta$ diminish (i.e. some of their elements are proven to be empty sets) from time $t_l$ to time $t_{l+1}$, then one applies the DSm hybrid rule in order to transfer the masses of empty sets to the non-empty sets (in the DSm classic rule's way) getting an updated basic belief assignment $m_{t_{l+1}|t_l}(.)$. Contrarily, if $\Theta$ and $D^\Theta$ increase (i.e. new elements arise in $\Theta$, and/or new elements in $D^\Theta$ are proven different from the empty set and as a consequence a basic belief assignment for them is required), then new masses (from the same or from the other sources of information) are needed to describe these new elements, and again one combines them using the DSm hybrid rule.

\clearpage
\newpage

%************************
\subsection{Example 3}
%************************

Let consider a fusion problem at time $t_l$ characterized by the frame $\Theta(t_l)\triangleq \{\theta_1,\theta_2\}$ and two independent sources of information providing the basic belief 
assignments $m_1(.)$ and $m_2(.)$ over $D^{\Theta(t_l)}$ and assume that at time $t_{l+1}$ a new hypothesis $\theta_3$ is introduced into the previous frame $\Theta(t_l)$ and a third 
source of evidence available at time $t_{l+1}$ provides its own basic belief assignment $m_3(.)$ over $D^{\Theta(t_{l+1})}$ where $$\Theta(t_{l+1})\triangleq \{ \Theta(t_l),\theta_3\} \equiv \{ \theta_1,\theta_2,\theta_3\}$$

\noindent
To solve such kind of dynamical fusion problems, we just use the classical DSm fusion rule  as follows:

\begin{itemize}
\item combine $m_1(.)$ and $m_2(.)$ at time $t_l$ using classical DSm fusion rule to get $m_{12}(.)= [m_1\oplus m_2](.)$ over $D^{\Theta(t_l)}$
\item because $D^{\Theta(t_l)} \subset D^{\Theta(t_{l+1})}$, $m_{12}(.)$ assigns the combined basic belief on a 
subset of $D^{\Theta(t_{l+1})}$, it is still directly possible to combine $m_{12}(.)$ with $m_3(.)$ at time $t_{l+1}$ by  the classical DSm fusion rule to get the final result $m_{123}(.)$ over $D^{\Theta(t_{l+1})}$ given by
$$m_{t_{l+1}}(.)\triangleq m_{123}(.)=  [m_{12}\oplus m_3](.) = [(m_1\oplus m_2)\oplus m_3](.) \equiv [m_1\oplus m_2\oplus m_3](.)$$
\item eventually apply DSm hybrid rule if some integrity constraints have to be taken into account in the model $\mathcal{M}$ of the problem
\end{itemize}

This method can be directly generalized to any number of sources of evidences and, in theory, to any structures/dimension of the frames $\Theta(t_{l})$, $\Theta(t_{l+1})$, ... In practice however, due to the huge number of elements of hyper-power sets,  the dimension of the frames $\Theta(t_{l})$, $\Theta(t_{l+1})$, $\dots$ must be not too large. This practical limitation depends on the computer resources available for the real-time processing. Specific suboptimal implementations of DSm rule will have to be developed to deal with fusion problems of large dimension. \\

It is also important to point out here that DSmT can easily deal, not only with dynamical fusion problems but with decentralized fusion problems as well working on non exhaustive frames. For example, let consider a set of two independent sources of information providing the basic belief 
assignments $m_1(.)$ and $m_2(.)$ over $D^{\Theta_{12}(t_l)=\{\theta_1,\theta_2\}}$ and another group of three independent sources of information providing the basic belief 
assignments $m_3(.)$, $m_4(.)$ and $m_5(.)$ over $D^{\Theta_{345}(t_l)=\{\theta_3,\theta_4,\theta_5,\theta_6\}}$, then it is still possible to combine all information in a decentralized manner as follows:

\begin{itemize}
\item combine $m_1(.)$ and $m_2(.)$ at time $t_l$ using classical DSm fusion rule to get $m_{12}(.)= [m_1\oplus m_2](.)$ over $D^{\Theta_{12}(t_l)}$.
\item combine $m_3(.)$, $m_4(.)$ and $m_5(.)$ at time $t_l$ using classical DSm fusion rule to get $m_{345}(.)= [m_3\oplus m_4\oplus m_5](.)$ over $D^{\Theta_{345}(t_l)}$.
\item consider now the global frame $\Theta(t_l)\triangleq \{\Theta_{12}(t_l) ,\Theta_{345}(t_l)\}$.
\item eventually apply DSm hybrid rule if some integrity constraints have to be taken into account in the model $\mathcal{M}$ of the problem.
\end{itemize}

Note that this {\it{static}} decentralized fusion can also be extended to decentralized dynamical fusion also by mixing two previous approaches.\\

One can even combine all five masses together by extending the vectors $m_i(.)$, $1\leq i \leq 5$, with null components for the new elements arisen from enlarging $\Theta$ to 
$\{\theta_1,\theta_2, \theta_3,\theta_4,\theta_5\}$ and correspondingly enlarging $D^\Theta$, and using the DSm hybrid rule for $k=5$. And more general combining the masses of any $k \geq 2$ sources.\\

We give now several simple numerical examples for such dynamical fusion problems involving non exclusive frames.

\subsubsection{Example 3.1}
%******************************

Let consider $\Theta(t_l)\triangleq \{\theta_1,\theta_2\}$ and the two following basic belief assignments available at time $t_l$:

$$m_1(\theta_1)=0.1 \qquad m_1(\theta_2)=0.2 \qquad m_1(\theta_1\cup\theta_2)=0.3\qquad m_1(\theta_1\cap\theta_2)=0.4$$
$$m_2(\theta_1)=0.5 \qquad m_2(\theta_2)=0.3 \qquad m_2(\theta_1\cup\theta_2)=0.1\qquad m_2(\theta_1\cap\theta_2)=0.1$$

The classical DSm rule of combination gives

$$m_{12}(\theta_1)=0.21 \qquad m_{12}(\theta_2)=0.17 \qquad m_{12}(\theta_1\cup\theta_2)=0.03\qquad m_{12}(\theta_1\cap\theta_2)=0.59$$

Now let consider at time $t_{l+1}$ the frame $\Theta(t_{l+1})\triangleq \{\theta_1,\theta_2,\theta_3\}$ and a third source of evidence with the following basic belief assignment

$$m_3(\theta_3)=0.4 \qquad  m_3(\theta_1\cap\theta_3)=0.3\qquad m_3(\theta_2\cup\theta_3)=0.3$$

\noindent
Then the final result of the fusion is obtained by combining $m_3(.)$ with $m_{12}(.)$ by the classical DSm rule of combination. One thus obtains:
\begin{multline*}
m_{123}(\theta_1\cap\theta_2\cap\theta_3)=0.464 \quad m_{123}(\theta_2\cap\theta_3)=0.068 \quad m_{123}(\theta_1\cap\theta_3)=0.156\quad m_{123}((\theta_1\cup \theta_2)\cap\theta_3)=0.012\\
m_{123}(\theta_1\cap\theta_2)=0.177  \quad m_{123}(\theta_1\cap(\theta_2\cup\theta_3))=0.063  \quad m_{123}(\theta_2)=0.051 \quad m_{123}((\theta_1\cap \theta_3)\cup\theta_2)=0.009
\end{multline*}

\subsubsection{Example 3.2}
%******************************

Let consider $\Theta(t_l)\triangleq \{\theta_1,\theta_2\}$ and the two previous following basic belief assignments $m_1(.)$ and $m_2(.)$ available at time $t_l$. The classical DSm fusion rule   
gives  gives as before

$$m_{12}(\theta_1)=0.21 \qquad m_{12}(\theta_2)=0.17 \qquad m_{12}(\theta_1\cup\theta_2)=0.03\qquad m_{12}(\theta_1\cap\theta_2)=0.59$$

Now let consider at time $t_{l+1}$ the frame $\Theta(t_{l+1})\triangleq \{\theta_1,\theta_2,\theta_3\}$ and the third source of evidence as in previous example with the basic belief assignment

$$m_3(\theta_3)=0.4 \qquad m_3(\theta_1\cap\theta_3)=0.3 \qquad m_3(\theta_2\cup\theta_3)=0.3 $$

\noindent
The final result of the fusion obtained by the classical DSm rule of combination corresponds to the result of the previous example, but suppose now one finds out that the integrity constraint $\theta_3= \emptyset$ holds, which implies also constraints $\theta_1\cap\theta_2\cap\theta_3=\emptyset$, $\theta_1\cap\theta_3=\emptyset$, $\theta_2\cap\theta_3=\emptyset$ and $(\theta_1\cup\theta_2)\cap\theta_3=\emptyset$. This is the DSm hybrid model $\mathcal{M}$ under consideration here. We then have to readjust the mass $m_{123}(.)$ of the previous example by the DSm hybrid rule and one finally gets

$$m_{\mathcal{M}}(\theta_1)=0.147$$ 
$$m_{\mathcal{M}}(\theta_2)=0.060+0.119=0.179$$
$$m_{\mathcal{M}}(\theta_1\cup\theta_2)=0 + 0 + 0.021=0.021$$
$$m_{\mathcal{M}}(\theta_1\cap\theta_2)=0.240+0.413=0.653$$

Therefore, when we restrain back $\theta_3= \emptyset$ and apply the DSm hybrid rule, we don't get back the same result (i.e. $m_{\mathcal{M}}(.)\neq m_{12}(.)$) because still remains some information from $m_3(.)$ on $\theta_1$, $\theta_2$, $\theta_1\cup\theta_2$, or $\theta_1\cap\theta_2$, i.e. $m_3(\theta_2)=0.3>0$.

\subsubsection{Example 3.3}
%******************************

Let consider $\Theta(t_l)\triangleq \{\theta_1,\theta_2\}$ and two previous following basic belief assignments $m_1(.)$ and $m_2(.)$ available at time $t_l$. The  classical DSm fusion rule gives as before
$$m_{12}(\theta_1)=0.21 \qquad m_{12}(\theta_2)=0.17 \qquad m_{12}(\theta_1\cup\theta_2)=0.03\qquad m_{12}(\theta_1\cap\theta_2)=0.59$$

Now let consider at time $t_{l+1}$ the frame $\Theta(t_{l+1})\triangleq \{\theta_1,\theta_2,\theta_3,\theta_4\}$ and another third source of evidence with the following basic belief assignment
$$m_3(\theta_3)=0.5 \qquad m_3(\theta_4)=0.3 \qquad m_3(\theta_3\cap\theta_4)=0.1 \qquad m_3(\theta_3\cup\theta_4)=0.1 $$
\noindent
Then, the DSm rule applied at time $t_{l+1}$ provides the following combined belief assignment
\begin{multline*}
m_{123}(\theta_1\cap\theta_3)=0.105 \quad m_{123}(\theta_1\cap\theta_4)=0.063 \quad m_{123}(\theta_1\cap(\theta_3\cup\theta_4))=0.021 \quad m_{123}(\theta_1\cap\theta_3\cap\theta_4)=0.021\\
m_{123}(\theta_2\cap\theta_3)=0.085 \quad m_{123}(\theta_2\cap\theta_4)=0.051 \quad m_{123}(\theta_2\cap(\theta_3\cup\theta_4))=0.017 \qquad m_{123}(\theta_2\cap\theta_3\cap\theta_4)=0.017 \\
m_{123}(\theta_3\cap(\theta_1\cup\theta_2))=0.015 \quad m_{123}(\theta_4\cap(\theta_1\cup\theta_2))=0.009 \quad 
m_{123}((\theta_1\cup\theta_2)\cap(\theta_3\cup\theta_4))=0.003\\ 
m_{123}((\theta_1\cup\theta_2)\cap(\theta_3\cap\theta_4))=0.003\quad  
m_{123}(\theta_1\cap\theta_2\cap\theta_3)=0.295 \quad
m_{123}(\theta_1\cap\theta_2\cap\theta_4)=0.177 \\
m_{123}((\theta_1\cap\theta_2)\cap(\theta_3\cup\theta_4))=0.059 \quad 
m_{123}(\theta_1\cap\theta_2\cap\theta_3\cap\theta_4)=0.059
\end{multline*}

Now suppose at time $t_{l+2}$ one finds out that $\theta_3=\theta_4=\emptyset$, then one applies the DSm hybrid rule after re-adjusting the combined belief mass $m_{123}(.)$ by cummulating the masses of  all empty sets. Using the DSm hybrid rule, one finally gets:

\begin{align*}
m_{t_{l+2}}(\theta_1)&=m_{123}(\theta_1)+\{m_{12}(\theta_1)m_3(\theta_3) + m_{12}(\theta_1)m_3(\theta_4) + m_{12}(\theta_1)m_3(\theta_3\cup\theta_4) + m_{12}(\theta_1)m_3(\theta_3\cap\theta_4)\}\\
&=0 + \{(0.21\times 0.5) + (0.21\times 0.3) + (0.21\times 0.1) + (0.21\times 0.1)\} = 0.21
\end{align*}
\begin{align*}
m_{t_{l+2}}(\theta_2)&=m_{123}(\theta_2)+\{m_{12}(\theta_2)m_3(\theta_3) + m_{12}(\theta_2)m_3(\theta_4) + m_{12}(\theta_2)m_3(\theta_3\cup\theta_4) + m_{12}(\theta_2)m_3(\theta_3\cap\theta_4)\}\\
&=0 + \{(0.17\times 0.5) + (0.17\times 0.3) + (0.17\times 0.1) + (0.17\times 0.1)\} = 0.17
\end{align*}
\begin{align*}
m_{t_{l+2}}(\theta_1\cup\theta_2)&=m_{123}(\theta_1\cup\theta_2)+\{m_{12}(\theta_1\cup\theta_2)m_3(\theta_3) + m_{12}(\theta_1\cup\theta_2)m_3(\theta_4) \\
&\qquad + m_{12}(\theta_1\cup\theta_2)m_3(\theta_3\cup\theta_4) + m_{12}(\theta_1\cup\theta_2)m_3(\theta_3\cap\theta_4)\}\\
&\qquad + \sum_{X_1, X_2 \in \{\theta_3, \theta_4, \theta_3\cup\theta_4, \theta_3\cap\theta_4\} } m_{12}(X_1)m_3(X_2)\\
&=0 + \{(0.03\times 0.5) + (0.03\times 0.3) + (0.03\times 0.1) +( 0.03\times 0.1)\} + \{0\} = 0.03
\end{align*}
\begin{align*}
m_{t_{l+2}}(\theta_1\cap\theta_2)&=m_{123}(\theta_1\cap\theta_2)+\{m_{12}(\theta_1\cap\theta_2)m_3(\theta_3) + m_{12}(\theta_1\cap\theta_2)m_3(\theta_4) \\
&\qquad + m_{12}(\theta_1\cap\theta_2)m_3(\theta_3\cup\theta_4) + m_{12}(\theta_1\cap\theta_2)m_3(\theta_3\cap\theta_4)\}\\
&=0 + \{(0.59\times 0.5) + (0.59\times 0.3) + (0.59\times 0.1) + (0.59\times 0.1)\} = 0.59
\end{align*}

\noindent
Thus we get the same result as for $m_{12}(.)$ at time $t_l$, which is normal.\\

{\bf{Remark:}}
note that if the third source of information don't assign non-null masses to $\theta_1$, or $\theta_2$ (or to their combinations using $\cup$ or $\cap$ operators), then one obtains the same result at time $t_{l+2}$ as at time $t_l$ as in this example 3.3, i.e. $m_{l+2}(.) = m_l(.)$, when imposing back $\theta_3=\theta_4=\emptyset$. But, if the third source of information assigns non-null masses to either $\theta_1$, or $\theta_2$, or to some of their combinations $\theta_1\cup \theta_2$ or $\theta_1\cap \theta_2$, then when one returns from 4 singletons to 2 singletons for $\Theta$, replacing $\theta_3=\theta_4=\emptyset$ and using the DSm hybrid rule, the fusion results at time $t_{l+2}$ is different from that at time $t_l$, and this is normal because some information/mass is left from the third source and is now fusioned with that of the previous sources (as in example 3.2 or in next example 3.4).\\

In general, let's suppose that the fusion of $k \geq 2$ masses provided by the sources $\mathcal{B}_1$, $\mathcal{B}_2$, ..., $\mathcal{B}_k$ has been done at time $t_l$ on $\Theta(t_l) = \{\theta_1, \theta_2, ..., \theta_n\}$.  At time $t_{l+1}$ new non-empty elements $\theta_{n+1}$, $\theta_{n+2}$, $\ldots$, $\theta_{n+m}$ appear, $m\geq 1$, thus $$\Theta(t_{l+1}) = \{\theta_1, \theta_2, ..., \theta_n,\theta_{n+1},\theta_{n+2}, \ldots, \theta_{n+m}\}$$
\noindent
and of course one or more sources (i.e. bodies of evidences) $\mathcal{B}_{k+1}$, $\ldots$, $\mathcal{B}_{k+l}$, where $l\geq 1$, appear to assign masses to these new elements.
\begin{itemize}
\item[a)] If all these new sources $\mathcal{B}_{k+1}$, $\ldots$, $\mathcal{B}_{k+l}$ assign null masses to all elements from $D^{\Theta(t_{l+1})}$ which contain in their structure/composition at least one of the singletons $\theta_1$, $\theta_2$, $\ldots$, $\theta_n$, then at time $t_{l+2}$ if one sets back the constraints that $\theta_{n+1}=\theta_{n+2}=\ldots = \theta_{n+m}=\emptyset$, then using the DSm hybrid rule, one obtains the same result as at time $t_l$, i.e. $m_{l+2}(.) = m_l(.)$.
\item[b)] Otherwise, the fusion at time $t_{l+2}$ will be different from the fusion at time $t_l$ because still remains some information/mass from sources $\mathcal{B}_{k+1}$, $\ldots$, $\mathcal{B}_{k+l}$ on singletons $\theta_1$, $\theta_2$, $\ldots$, $\theta_n$ or on some elements from $D^{\Theta(t_l)}$ which contain at least one of such singletons, information/mass which fusions with the previous sources.
\end{itemize}

\subsubsection{Example 3.4}
%******************************

Let consider $\Theta(t_l)\triangleq \{\theta_1,\theta_2\}$ and the two following basic belief assignments available at time $t_l$:
$$m_1(\theta_1)=0.6 \qquad m_1(\theta_2)=0.4 \qquad\text{and}\qquad m_2(\theta_1)=0.7 \qquad m_2(\theta_2)=0.3 $$
\noindent The classical DSm rule of combination gives $m_{12}(\theta_1)=0.42$, $m_{12}(\theta_2)=0.12$ and  $m_{12}(\theta_1\cap\theta_2)=0.46$. Now let consider at time $t_{l+1}$ the frame $\Theta(t_{l+1})\triangleq \{\theta_1,\theta_2,\theta_3\}$ and a third source of evidence with the following basic belief assignment $m_3(\theta_1)=0.5$, $m_3(\theta_2)=0.2$ and $m_3(\theta_3)=0.3$. Then the final result obtained from the classical DSm rule of combination is still as before
\begin{multline*}
m_{123}(\theta_1)=0.210 \quad m_{123}(\theta_2)=0.024 \quad m_{123}(\theta_1\cap\theta_2)=0.466 \quad m_{123}(\theta_1\cap\theta_3)=0.126\\
m_{123}(\theta_2\cap\theta_3)=0.036 \quad
m_{123}(\theta_1\cap\theta_2\cap\theta_3)=0.138
\end{multline*}

Suppose now one finds out that the integrity constraint $\theta_1\cap\theta_3= \emptyset$ which also implies $\theta_1\cap\theta_2\cap\theta_3=\emptyset$. This is the DSm hybrid model $\mathcal{M}$ under consideration. By applying the DSm hybrid fusion rule, one forces 
$m_{\mathcal{M}}(\theta_1\cap\theta_3)=0$ and $m_{\mathcal{M}}(\theta_1\cap\theta_2\cap\theta_3)=0$ and we transfer $m_{123}(\theta_1\cap\theta_2\cap\theta_3)=0.138$ towards $m_{\mathcal{M}}((\theta_1\cap\theta_2)\cup\theta_3)$ and the mass $m_{123}(\theta_1\cap\theta_3)=0.126$ has to be transferred towards
$m_{\mathcal{M}}(\theta_1\cup\theta_3)$.
One then gets finally 
\begin{multline*}
m_{\mathcal{M}}(\theta_1)=0.210 \quad m_{\mathcal{M}}(\theta_2)=0.024 \quad m_{\mathcal{M}}(\theta_1\cap\theta_2)=0.466 \quad
m_{\mathcal{M}}(\theta_2\cap\theta_3)=0.036 \\
m_{\mathcal{M}}((\theta_1\cap\theta_2)\cup\theta_3)=0.138 \quad
m_{\mathcal{M}}(\theta_1\cup\theta_3)=0.126
\end{multline*}

\subsubsection{Example 3.5}
%******************************

Let consider $\Theta(t_l)\triangleq \{\theta_1,\theta_2\}$ and the two previous  basic belief assignments available at time $t_l$ as in previous example, i.e.
$$m_1(\theta_1)=0.6 \qquad m_1(\theta_2)=0.4 \qquad\text{and}\qquad m_2(\theta_1)=0.7 \qquad m_2(\theta_2)=0.3 $$
\noindent The classical DSm rule of combination gives
$$m_{12}(\theta_1)=0.42 \qquad m_{12}(\theta_2)=0.12 \quad m_{12}(\theta_1\cap\theta_2)=0.46$$
\noindent Now let consider at time $t_{l+1}$ the frame $\Theta(t_{l+1})\triangleq \{\theta_1,\theta_2,\theta_3\}$ and a third source of evidence with the following basic belief assignment
$$m_3(\theta_1)=0.5 \qquad m_3(\theta_2)=0.2 \qquad m_3(\theta_3)=0.3 $$

\noindent
Then the final result of the fusion is obtained by combining $m_3(.)$ with $m_{12}(.)$ by the classical DSm rule of combination. One thus obtains now
\begin{multline*}
m_{123}(\theta_1)=0.210 \quad m_{123}(\theta_2)=0.024 \quad m_{123}(\theta_1\cap\theta_2)=0.466 \quad m_{123}(\theta_1\cap\theta_3)=0.126\\  m_{123}(\theta_2\cap\theta_3)=0.036 \quad m_{123}(\theta_1\cap\theta_2\cap\theta_3)=0.138
\end{multline*}

\noindent
But suppose one finds out that the integrity constraint is now $\theta_3=\emptyset$ which implies necessarily also $\theta_1\cap\theta_3=\theta_2\cap\theta_3=\theta_1\cap\theta_2\cap\theta_3\equiv\emptyset$ and $(\theta_1\cup\theta_2)\cap\theta_3=\emptyset$ 
(this is our new DSm hybrid model $\mathcal{M}$ under consideration in this example). By applying the DSm hybrid fusion rule, one gets finally the non-null masses 
$$m_{\mathcal{M}}(\theta_1)=0.336 \qquad m_{\mathcal{M}}(\theta_2)=0.060\qquad 
m_{\mathcal{M}}(\theta_1\cap\theta_2)=0.604$$

\subsubsection{Example 3.6}
%******************************

Let consider $\Theta(t_l)\triangleq \{\theta_1,\theta_2,\theta_3,\theta_4\}$ and the following basic belief assignments available at time $t_l$ :
$$\begin{cases}
m_1(\theta_1)=0.5 \qquad m_1(\theta_2)=0.4 \qquad m_1(\theta_1\cap\theta_2)=0.1\\
m_2(\theta_1)=0.3 \qquad m_2(\theta_2)=0.2 \qquad m_2(\theta_1\cap\theta_3)=0.1\qquad m_2(\theta_4)=0.4
\end{cases}
$$
\noindent The classical DSm rule of combination gives
$$m_{12}(\theta_1)=0.15 \qquad m_{12}(\theta_2)=0.08 \qquad m_{12}(\theta_1\cap\theta_2)=0.27 \qquad m_{12}(\theta_1\cap\theta_3)=0.05\qquad m_{12}(\theta_1\cap\theta_4)=0.20$$ $$m_{12}(\theta_2\cap\theta_4)=0.16 \qquad m_{12}(\theta_1\cap\theta_2\cap\theta_3)=0.05\qquad m_{12}(\theta_1\cap\theta_2\cap\theta_4)=0.04$$
\noindent
Now assume that at time $t_{l+1}$ one finds out that $\theta_1\cap\theta_2\overset{\mathcal{M}}{\equiv}\theta_1\cap\theta_3\overset{\mathcal{M}}{\equiv}\emptyset$. 
Using the DSm hybrid rule, one gets:
$$
\begin{cases}
m_{\mathcal{M}}(\theta_1\cap\theta_2)=m_{\mathcal{M}}(\theta_1\cap\theta_3)=m_{\mathcal{M}}(\theta_1\cap\theta_2\cap\theta_3)=m_{\mathcal{M}}(\theta_1\cap\theta_2\cap\theta_4)=0\\
m_{\mathcal{M}}(\theta_1)=m_{12}(\theta_1)+m_2(\theta_1)m_1(\theta_1\cap\theta_2)+m_1(\theta_1)m_2(\theta_1\cap\theta_3)=0.15+0.03+0.05=0.23\\
m_{\mathcal{M}}(\theta_2)=m_{12}(\theta_2)+m_2(\theta_2)m_1(\theta_1\cap\theta_2)+m_1(\theta_2)m_2(\theta_1\cap\theta_3)=0.08+0.02+0.04=0.14\\
m_{\mathcal{M}}(\theta_4)=m_{12}(\theta_4)+m_1(\theta_1\cap\theta_2)m_2(\theta_4)=0+0.04=0.04\\
m_{\mathcal{M}}(\theta_1\cap\theta_4)=m_{12}(\theta_1\cap\theta_4)=0.20\\
m_{\mathcal{M}}(\theta_2\cap\theta_4)=m_{12}(\theta_2\cap\theta_4)=0.16\\
m_{\mathcal{M}}(\theta_1\cup\theta_2)=m_{12}(\theta_1\cup\theta_2)+m_1(\theta_1)m_2(\theta_2)+m_2(\theta_1)m_1(\theta_2)+m_1(\theta_1\cap\theta_2)m_2(\theta_1\cap\theta_2)=
0.22\\
m_{\mathcal{M}}(\theta_1\cup\theta_2\cup\theta_3)=m_{12}(\theta_1\cup\theta_2\cup\theta_3)+m_1(\theta_1\cap\theta_2)m_2(\theta_1\cap\theta_3)+m_2(\theta_1\cap\theta_2)m_1(\theta_1\cap\theta_3)\\
\qquad\qquad\qquad\qquad\quad+m_1(\theta_1\cap\theta_2\cap\theta_3)m_2(\theta_1\cap\theta_2\cap\theta_3)=0.01
\end{cases}
$$

\subsubsection{Example 3.7}
%******************************

Let consider $\Theta(t_l)\triangleq \{\theta_1,\theta_2,\theta_3,\theta_4\}$ and the following basic belief assignments available at time $t_l$ :
$$\begin{cases}
m_1(\theta_1)=0.2 \qquad m_1(\theta_2)=0.4 \qquad m_1(\theta_1\cap\theta_2)=0.1 \qquad m_1(\theta_1\cap\theta_3)=0.2\qquad m_1(\theta_4)=0.1\\
m_2(\theta_1)=0.1 \qquad m_2(\theta_2)=0.3 \qquad m_2(\theta_1\cap\theta_2)=0.2 \qquad m_2(\theta_1\cap\theta_3)=0.1\qquad m_2(\theta_4)=0.3
\end{cases}
$$
\noindent The classical DSm rule of combination gives
$$m_{12}(\theta_1)=0.02 \qquad m_{12}(\theta_2)=0.12 \qquad m_{12}(\theta_1\cap\theta_2)=0.28 \qquad m_{12}(\theta_1\cap\theta_3)=0.06\qquad m_{12}(\theta_4)= 0.03$$
$$m_{12}(\theta_1\cap\theta_4)=0.07\qquad m_{12}(\theta_2\cap\theta_4)=0.15 \qquad m_{12}(\theta_1\cap\theta_2\cap\theta_3)=0.15$$
$$m_{12}(\theta_1\cap\theta_2\cap\theta_4)=0.05\qquad m_{12}(\theta_1\cap\theta_3\cap\theta_4)=0.07$$
\noindent
Now assume that at time $t_{l+1}$ one finds out that $\theta_1\cap\theta_2\overset{\mathcal{M}}{\equiv}\theta_1\cap\theta_3\overset{\mathcal{M}}{\equiv}\emptyset$. 
Using the DSm hybrid rule, one gets:
$$
\begin{cases}
m_{\mathcal{M}}(\theta_1\cap\theta_2)=m_{\mathcal{M}}(\theta_1\cap\theta_3)=m_{\mathcal{M}}(\theta_1\cap\theta_2\cap\theta_3)=m_{\mathcal{M}}(\theta_1\cap\theta_2\cap\theta_4)=0\\
m_{\mathcal{M}}(\theta_1)=m_{12}(\theta_1)+m_1(\theta_1)m_2(\theta_1\cap\theta_2)+m_2(\theta_1)m_1(\theta_1\cap\theta_2)+m_1(\theta_1)m_2(\theta_1\cap\theta_3)+m_2(\theta_1)m_1(\theta_1\cap\theta_3)=
0.11\\
m_{\mathcal{M}}(\theta_2)=m_{12}(\theta_2)+m_1(\theta_2)m_2(\theta_1\cap\theta_2)+m_2(\theta_2)m_1(\theta_1\cap\theta_2)+m_1(\theta_2)m_2(\theta_1\cap\theta_3)+m_2(\theta_2)m_1(\theta_1\cap\theta_3)=
0.33\\
m_{\mathcal{M}}(\theta_4)=m_{12}(\theta_4)+m_1(\theta_4)m_2(\theta_1\cap\theta_2)+m_2(\theta_4)m_1(\theta_1\cap\theta_2)+m_1(\theta_4)m_2(\theta_1\cap\theta_3)+m_2(\theta_4)m_1(\theta_1\cap\theta_3)=
0.15\\
m_{\mathcal{M}}(\theta_1\cap\theta_4)=m_{12}(\theta_1\cap\theta_4)=0.07\\
m_{\mathcal{M}}(\theta_2\cap\theta_4)=m_{12}(\theta_2\cap\theta_4)=0.15\\
m_{\mathcal{M}}(\theta_1\cup\theta_2)=m_{12}(\theta_1\cup\theta_2)+m_1(\theta_1\cap\theta_2)m_2(\theta_1\cap\theta_2)+m_1(\theta_1)m_2(\theta_2)+m_2(\theta_1)m_1(\theta_2)=
0.12\\
m_{\mathcal{M}}(\theta_1\cup\theta_3)=m_{12}(\theta_1\cup\theta_3)+m_1(\theta_1\cap\theta_3)m_2(\theta_1\cap\theta_3)=
0.02\\
m_{\mathcal{M}}(\theta_1\cup\theta_2\cup\theta_3)=m_{12}(\theta_1\cup\theta_2\cup\theta_3)+m_1(\theta_1\cap\theta_2)m_2(\theta_1\cap\theta_3)+m_2(\theta_1\cap\theta_2)m_1(\theta_1\cap\theta_3)=
0.05
\end{cases}
$$

%----------------------------------------------------------------------------------------------------------------------------
\section{Bayesian mixture of DSm hybrid models}
%----------------------------------------------------------------------------------------------------------------------------

In the preceding, one has first shown how to combine generalized basic belief assignments provided by $k\geq 2$ independent sources of information with the general DSm hybrid rule of combination which deals with all possible kinds of constraints introduced by the hybrid model of the problem. This approach implicitly assumes that one knows/trusts with certainty that the model $\mathcal{M}$ (usually a DSm hybrid model) of the problem is valid and corresponds to the true model. In some complex fusion problems however (static or dynamic ones), one may have some doubts about the validity of the model $\mathcal{M}$ on which is based the fusion because of the nature and evolution of elements of the frame $\Theta$. In such situations, we propose to consider a set of exclusive and exhaustive models $\{\mathcal{M}_1,\mathcal{M}_2,\ldots,\mathcal{M}_K\}$ with some probabilities $\{P\{\mathcal{M}_1\},P\{\mathcal{M}_2\},\ldots,P\{\mathcal{M}_K\}\}$. We don't go here deeper on the justification/acquisition of such probabilities because this is highly dependent on the nature of the fusion problem under consideration. We just assume here that such probabilities are available at any given time $t_l$ when the fusion has to be done. We propose then to use the Bayesian mixture of combined masses $m_{\mathcal{M}_i(\Theta)}(.)$ $i=1,\ldots,K$ to obtain the final result :
\begin{equation}
\forall A\in D^\Theta, \qquad m_{\mathcal{M}_1,\ldots,\mathcal{M}_K}(A)=\sum_{i=1}^K  P\{\mathcal{M}_i\}  m_{\mathcal{M}_i(\Theta)}(A)
\end{equation}

\section{ Conclusion}
%--------------------------------------------------------------------------------------------

In this paper we have extended the DSmT and the classical DSm rule of combination to the case of  hybrid models for the frame of discernment involved in many complex fusion problems. The free-DSm model (which assumes that none of the elements of the frame is refinable) can be interpreted as the opposite of the Shafer's model (which assumes that all elements of the frame are truly exclusive) on which is based the mathematical theory of evidence (Dempster-Shafer Theory - DST). Between these two extreme models, there exists actually many possible hybrid models for the frames of discernment depending on the real intrinsic nature of elements of the fusion problem under consideration. For real problems, some elements of the frame of discernment can appear to be truly exclusive whereas some others cannot be considered as fully discernable or refinable. This present research work proposes a new  DSm hybrid rule of combination for hybrid-models based on the DSmT. The DSm hybrid rule works in any model and is involved in calculation of mass fusion of any number of sources of information, no matter how big is the
conflict/paradoxism of sources, and on any frame of discernment (exhaustive
or non-exhaustive, with elements which may be exclusive or non-exclusive
or both).  This is an important rule since does not require the calculation
of weighting factors, neither normalization as other rules do, and the
transfer of empty-sets' masses to non-empty sets masses is naturally done
following the DSm network architecture which is derived from the DSm classic
rule. DSmT together with DSm hybrid rule appears from now on to be a new alternative to classical approaches and to existing combination rules and to be very promising for the development of future complex (uncertain/incomplete/paradoxical/dynamical)  information fusion systems.

\end{document}